\documentclass[trsc,nonblindrev]{informs3} 

\OneAndAHalfSpacedXI 

\usepackage{hyperref}

\usepackage{etoolbox}
\makeatletter
\patchcmd{\maketitle}
 {\def\@makefnmark}
 {\def\@makefnmark{}\def\useless@macro}
 {}{}
\makeatother

\usepackage{subfig}

\usepackage{array,multirow,graphicx}

\usepackage{adjustbox}

\usepackage{natbib}
 \bibpunct[, ]{(}{)}{,}{a}{}{,}%
 \def\bibfont{\small}%
\TheoremsNumberedThrough     

\EquationsNumberedThrough    

\begin{document}


 \RUNAUTHOR{Papadopoulos et al.} 

\RUNTITLE{Personalized Pareto-Improving Pricing Schemes with Truthfulness Guarantees}

\TITLE{Personalized Pareto-Improving Pricing-and-Routing Schemes for Near-Optimum Freight Routing: An Alternative Approach to Congestion Pricing\thanks{Declarations of interest: none}}

\ARTICLEAUTHORS{%
\AUTHOR{Aristotelis-Angelos Papadopoulos}
\AFF{Ming Hsieh Department of Electrical and Computer Engineering, \\ University of Southern California, 
   Los Angeles, CA 90089 USA, \EMAIL{aristotp@usc.edu (corresponding author)} \URL{}}
\AUTHOR{Ioannis Kordonis}
\AFF{CentraleSup\'elec, Avenue de la Boulaie, 35576 Cesson-S\'evign\'e, France, \EMAIL{jkordonis1920@yahoo.com} \URL{}}
\AUTHOR{Maged M. Dessouky}
\AFF{Daniel J. Epstein Department of Industrial and Systems Engineering, \\ University of Southern California, Los Angeles, CA 90089 USA, \EMAIL{maged@usc.edu} \URL{}}
\AUTHOR{Petros A. Ioannou}
\AFF{Ming Hsieh Department of Electrical and Computer Engineering, \\ University of Southern California, 
   Los Angeles, CA 90089 USA, \EMAIL{ioannou@usc.edu} \URL{}}
} 

\ABSTRACT{%
Traffic congestion constitutes a major problem in urban areas. Trucks contribute to congestion and have a negative impact on the environment due to their size, slower dynamics and higher fuel consumption. The individual routing decisions made by truck drivers do not lead to system optimum operations and contribute to traffic imbalances especially in places where the volume of trucks is relatively high. In this paper, we design a coordination mechanism for truck drivers that uses pricing-and-routing schemes that can help alleviate traffic congestion in a general transportation network. We consider the user heterogeneity in Value-Of-Time (VOT) by adopting a multi-class model with stochastic Origin-Destination (OD) demands for the truck drivers. The main characteristic of the mechanism is that the coordinator asks the truck drivers to declare their desired OD pair and pick their individual VOT from a set of $N$ available options, and guarantees that the resulting pricing-and-routing scheme is Pareto-improving, i.e. every truck driver will be better-off compared to the User Equilibrium (UE) and that every truck driver will have an incentive to truthfully declare his/her VOT, while leading to a revenue-neutral (budget balanced) on average mechanism. This approach enables us to design personalized (VOT-based) pricing-and-routing schemes. We show that the Optimum Pricing Scheme (OPS) can be calculated by solving a nonconvex optimization problem. To improve computational efficiency, we propose an Approximately Optimum Pricing Scheme (AOPS) and prove that it satisfies the aforementioned properties. Both pricing-and-routing schemes are compared to the Congestion Pricing with Uniform Revenue Refunding (CPURR) scheme through extensive simulation experiments where it is shown that OPS and AOPS achieve a much lower expected total travel time and expected total monetary cost for the users compared to the CPURR scheme, without negatively affecting the rest of the network. These results demonstrate the efficiency of personalized (VOT-based) pricing-and-routing schemes.  
}%



\maketitle

%

\section{Introduction.}
Measuring the contribution to the United States (U.S.) economy as the share of all expenditures in transportation-related final goods and services, the transportation sector contributed $8.9\%$ to U.S. Gross Domestic Product (GDP) \citep{statistics1} while in the European Union (EU) it accounts for almost $5\%$ of the GDP \citep{statistics3}. In EU, road transport has the largest share of EU freight transport accounting
for 76.7\% of the total inland freight transport \citep{statistics4}. \cite{statistics2} found that the trucking industry experienced nearly 1.2 billion hours of delay on the National Highway System (NHS) of the U.S. as a result of traffic congestion making the operational costs incurred by the trucking industry due to traffic congestion to be $\$74.5$ billion per year. These statistics demonstrate that an optimized routing system is essential and could significantly contribute to the global economy. 

Drivers usually make their routing decisions using GPS routing apps in an effort to minimize their individual travel time or cost objective. This phenomenon is known as User Equilbrium (UE) or the first Wardrop Principle \citep{wardrop}. However, it is known that UE deviates from an optimized road usage \citep{beckmann, pigou} and it is a sub-optimal behavior compared to the socially optimum policy that could be achieved through a centrally coordinated system \citep{Youn_2008}. Recent studies \citep{monnot2017bad, Zhang_2016} estimated the Price Of Anarchy (POA) \citep{POA}, i.e. the inefficiency between a selfish routing strategy and a system optimum policy in realistic transportation networks using real traffic data, demonstrating the necessity for its reduction. Based on the idea of Connected Automated Vehicles (CAVs) \citep{malikopoulos, zhang12}, \cite{Zhang_2016} proposed to reduce the POA by recommending to all drivers socially optimum routes. However, such a strategy would raise several fairness and equity issues since in a System Optimum (SO) solution, some drivers may benefit while some others may be harmed compared to the UE.  

One of the most common techniques addressing the problem of the inefficiency between the UE and the SO solutions is congestion pricing \citep{cong_pric_2020, vickrey, beckmann, pigou} where each driver is assigned a fee corresponding to the additional cost his/her presence causes to the network. Several other works have studied congestion pricing under user heterogeneity in VOT, e.g. \citep{congest_pric_heter, congest_pric_heter222}, the problem of management of the revenue collected from the application of congestion pricing \citep{RePEc:eee:transb:v:44:y::i:8-9:p:972-982, revenues} and the impact of congestion pricing schemes on emissions of freight transport \citep{emissions}. London \citep{london_cong_pric}, Stockholm \citep{ELIASSON2009240}, Singapore and Milan \citep{LEHE2019200} are some of the cities that have already introduced congestion pricing, while recent studies \citep{cong_pric_2020_2, cong_pric_rome} also explore the benefits from applying congestion pricing to more major cities. Recently, there is also a growing research interest for studies related to pricing schemes in the presence of autonomous vehicles \citep{lazar2019optimal,mehr2019pricing, congestion3, congestion2}. 

Another well studied set of strategies addressing the problem of the inefficiency between an equilibrium flow pattern and the SO are the applications of Tradable Credit Schemes (TCS) \citep{itsc_2019_credits, credits1} or tradable travel permits \citep{WADA201394,permits_2020_implem, permits2} among the drivers of the network. In this case, a central coordinator initially distributes a certain number of credits (or permits) to all eligible drivers and free credit (or permit) trading is allowed among travelers. \cite{credits2} and \cite{credits3} studied the application of TCS under user heterogeneity in VOT, while \cite{permits1} studied OD-based travel permits in the presence of heterogeneous users. Recently, \cite{credits4} studied a Cyclic Tradable Credit Scheme (CTCS), where the credits never expire but circulate within the system, and derived a sufficient condition for the existence of a Pareto-improving CTCS in a general network. For a more comprehensive review of credit- and permit-based schemes, we refer the interested reader to \cite{review_permits}.

In this paper, we address the problem of the inefficiency of an equilibrium flow pattern by studying pricing schemes under a centrally coordinated freight routing system that can alleviate traffic congestion and drive the network as close as possible to a SO solution. We focus our study on pricing-and-routing schemes that can be specifically applied on trucks. Given that truck drivers routinely use varying routes for the same journey depending on the traffic conditions \citep{giannis} and the fact that their travel time is already a commodity, make trucks form an ideal candidate subclass of vehicles for coordinated routing. To this end, we consider a non-atomic game theoretic model whose users are the truck drivers and their demand is assumed to be stochastic. In the case where the planning horizon is split into discrete non-overlapping time intervals and the drivers choose both their OD pair as well as their desired departure time interval, \cite{papadopoulos2019incentives, aris_ecc} derived sufficient conditions for the existence of revenue-neutral (budget balanced) and Pareto-improving pricing schemes that can additionally provide individual incentives to the drivers to truthfully declare their desired departure time. In this work, we take into account the user heterogeneity in the VOT. For the single OD case, using a bottleneck model \citep{vickrey} and assuming two classes of users with distinct VOT, \cite{incentives_2020_tr_part_b} explored the possibility of adopting the instrument of incentives to shift commuters’ departure times in
a single morning bottleneck situation. For the fixed demand case, \cite{RePEc:eee:transb:v:44:y::i:8-9:p:972-982} derived sufficient conditions for the existence of Pareto-improving and revenue-neutral pricing schemes. However, since they could not find a way to identify the VOT of each user, they proposed class-anonymous pricing schemes based on the idea of Congestion Pricing with Uniform Revenue Refunding (CPURR).

\citep{vot_geroliminis} and \citep{ZHENG2016133} argued that VOT-based pricing schemes can increase the feasibility of implementation since they take into account the vulnerable user groups. However, most of the existing literature, e.g. \citep{e430d1dc1f7a428d83c3c49d1a22bcfc, RePEc:eee:transe:v:54:y:2013:i:c:p:1-13}, makes assumptions about the distribution that the VOT of the drivers might follow and to the best of our knowledge, no self-reporting scheme where the users directly report their VOT to a central authority has been previously proposed. Note that under such a scheme, it would be important to provide incentives to the users to truthfully report their VOT in order to avoid the exploitability of the mechanism. This is mainly because many users would be willing to declare a high VOT in order to be assigned to the fastest possible route. In this work, we design a coordination mechanism for the truck drivers where the central coordinator asks the users to declare their desired OD pair and additionally pick their VOT from a set of $N$ available options. Under this structure, we prove the existence of revenue-neutral and Pareto-improving pricing schemes in the sense that the expected cost of the users at the time they make their decision is less than or equal to their corresponding cost at the UE, that can additionally provide incentives to the drivers to truthfully declare their VOT. This additional information enables us to design personalized (VOT-based) pricing-and-routing schemes. More specifically, we propose an Optimum Pricing Scheme (OPS) that can be calculated by solving a nonconvex optimization problem. To reduce the computational time needed to calculate the OPS, we propose a second pricing-and-routing scheme called Approximately Optimum Pricing Scheme (AOPS) and we prove that it satisfies the desired properties. The simulation experiments demonstrate that both OPS and AOPS provide a much lower expected total travel time and expected total monetary cost to the users compared to the CPURR scheme, while concurrently approaching the SO solution. 

The rest of the paper is organized as follows. In Section 2, we present the model used and we formulate the User Equilibrium (UE) and the System Optimum (SO) problems. In Section 3, we present the Optimum Pricing Scheme (OPS) and the Approximately Optimum Pricing Scheme (AOPS) and we additionally formulate the CPURR scheme in the form of an optimization problem with complementarity constraints. In Section 4, the simulation results of our approach are provided while in Section 5, we present the Conclusion of this work. 

\section{Problem Formulation.}
Let $G=(V,L)$ denote a transportation network, where $V$ is the set of nodes and $L$ is the set of links in the network. Let $C_{lT}(X_{lp}, X_{lT}(\alpha))$ be a known nonlinear function representing the travel time of a truck driver traversing road segment $l$ when there exist $X_{lp}$ passenger vehicles and $X_{lT}(\alpha)$ trucks on it, where $\alpha$ is a set of variables defined as follows:
\begin{equation}
\begin{aligned}
\alpha = \{\alpha_{w,r}^{j}: w=1,\dots,N, j=1,\dots,v, r\in R_j\}
\end{aligned}
\end{equation}
where $j$ is the index corresponding to a specific Origin-Destination (OD) pair, $w$ is the index corresponding to a class of users with VOT $s_w$, $r \in R_j$ denotes a specific route among the set of available routes $R_j$ connecting OD pair $j$, $N$ is the number of distinct classes of users and $v$ is the number of OD pairs in the network. Therefore, $\alpha_{w,r}^{j}$ expresses the proportion of truck drivers belonging to class $w$ with a desired OD pair $j$ who choose route $r$ for their trip. Additionally, we assume that the OD demand of the truck drivers is stochastic and follows a probability distribution with finite support. Let $d_{j,w}$ be random variables denoting the number of truck drivers belonging to the class $w$ with desired OD pair $j$ and let $d_{j,w}^{c}$ be their corresponding values during the demand realization $c$. Then, the number of trucks traversing the road segment $l$ is given by:
\begin{equation}
\label{proti}
    X_{lT}(\alpha) = \sum_{j=1}^v \sum_{w=1}^N \sum_{r \in R_j: l \in r} d_{j,w}^{c} \alpha_{w,r}^j
\end{equation}
where on the left side of (\ref{proti}), we omitted the index $c$ to simplify the notation.

We consider a model with a continuum of users. We assume that truck drivers know the number of passenger vehicles at each road segment of the transportation network. In case they do not possess this information and they make their routing decisions based on the probability distribution of $X_{lp}$, we can write the travel time of each road segment $l$ as:
\begin{equation}
\label{double_expect}
    E\Bigg[\sum_{l \in r} C_{lT}(X_{lp}, X_{lT}(\alpha))\Bigg] = E\Bigg[\sum_{l \in r} E\Big[C_{lT}(X_{lp}, X_{lT}(\alpha))|X_{lT}\Big]\Bigg] 
\end{equation}
which is a function of $X_{lT}$ only, leading to a similar analysis. Note that in (\ref{double_expect}), the inner expectation is taken with respect to $X_{lp}$ and the outer expectation is taken with respect to $X_{lT}$. For the rest of the analysis, we assume that the number of passenger vehicles in the network $X_{lp}$ is a deterministic quantity and is known by the truck drivers. Additionally, we assume that the truck drivers know the probability distribution of the OD demand for the rest of the truck drivers but not the exact realization of the demand. This is a symmetric information model since each truck driver has the same amount of information. A similar model was also used in \cite{giannis}. In this paper, we extend this model by considering user heterogeneity in VOT of the truck drivers.

The above formulation is used in subsequent sections to study the User Equilibrium (UE) and System Optimum (SO) flow patterns as well as the existence of Pareto-improving pricing-and-routing schemes.

\subsection{User Equilibrium (UE).}
In the absence of pricing schemes, the drivers are trying to minimize their own individual travel time, e.g. through the usage of GPS routing apps. This behavior drives the network to a state called User Equilibrium \citep{wardrop} where no driver has an incentive to unilaterally change his/her routing decision since he/she is not going to benefit from such a change.

In a transportation network with heterogeneous users, the equilibrium conditions can be either calculated in time units or in cost units \citep{RePEc:eee:transb:v:44:y::i:8-9:p:972-982}. Since the equilibrium conditions expressed in cost units can be obtained by multiplying each class' travel time by its corresponding VOT, we formulate the UE problem in time units without any loss of generality. Therefore, let $F_{j,r}^{w}(\alpha)$ be the expected travel time of a truck driver with VOT belonging to the class $w$, travelling in OD pair $j$ and following route $r$. Then, $F_{j,r}^{w}(\alpha)$ is given by:
\begin{equation}
\label{help1}
\begin{aligned}
F_{j,r}^{w}(\alpha) = E\Bigg[\sum_{l \in r} C_{lT}(X_{lp}, X_{lT}(\alpha))\Bigg]
\end{aligned}
\end{equation}
where $X_{lT}(\alpha)$ is given by (\ref{proti}). Note that in a UE solution, it holds that $F_{j,r}^{w}(\alpha) = F_{j,r}^{w^{\prime}}(\alpha), \forall w \neq w^{\prime}$, i.e. the equilibrium travel time is identical for all user classes between the same OD pair. Additionally, in an equilibrium condition, it holds that:
\begin{equation}
\label{equil_cond}
\begin{aligned}
\alpha^{j}_{w,r} > 0 \implies F_{j,r}^{w}(\alpha) \leq F_{j,r'}^{w}(\alpha), \forall  r' \neq r
\end{aligned}
\end{equation}
where $r,r' \in R_j$. Inequality (\ref{equil_cond}) states that in an equilibrium condition, drivers are choosing the route $r$ that minimizes their individual expected travel time. 

It has been shown by \cite{giannis} that there are possibly many non-equivalent UE solutions. In this work, we calculate an equilibrium solution by solving an optimization problem with complementarity constraints \citep{pang_complementarity} which is a nonconvex optimization problem. Before formulating the problem, let us first define the expected total travel time of the truck drivers in the network as:
\begin{equation}
\label{trav_time_tr}
\begin{aligned}
E[T_{tr}(\alpha)] = E\Bigg[\sum_{l=1}^{m} X_{lT}(\alpha) C_{lT}(X_{lp}, X_{lT}(\alpha))\Bigg]
\end{aligned}
\end{equation}
where $m$ is the number of road segments in the transportation network and $X_{lT}(\alpha)$ is given by (\ref{proti}). Under the assumption that the demand of the truck drivers follows a probability distribution with finite support, we can define the expected total monetary cost of the truck drivers in the network as:
\begin{equation}
\label{mon_cost}
\begin{aligned}
E[T_{tr}^{mon}(\alpha)] = \sum_c \sum_{j=1}^{v} \sum_{w=1}^{N} \sum_{r \in R_j} p_c d_{j,w}^{c} \alpha_{j,w,r}^{UE} s_w J_{c,j,w,r}^{UE}(\alpha)
\end{aligned}
\end{equation}
where $c$ and $p_c$ correspond to a specific realization of the demand $d_{j,w}^{c}$ and its associated probability, respectively. Note that since we assumed that the demand of the truck drivers follows a probability distribution with finite support, as it is common in probability theory \citep{10.5555/1207291}, we use the term `demand realization' to describe the observed value that the demand of the truck drivers takes. Moreover, $\alpha_{j,w,r}^{UE}$ is the proportion of truck drivers belonging to class $w$ with a desired OD pair $j$ who follow route $r$ at the UE, $s_w$ is the VOT of the class $w$ and $J_{c,j,w,r}^{UE}$ is the travel time of a truck driver with an OD pair $j$ who follows route $r$ during the demand realization $c$ at the UE and is given by the following equation:
\begin{equation}
\label{J_UE}
    J_{c,j,w,r}^{UE}(\alpha) = \sum_{l \in r} C_{lT}(X_{lp}, X_{lT}(\alpha))
\end{equation}
where $X_{lT}(\alpha)$ is given by (\ref{proti}). Using (\ref{J_UE}), we can rewrite (\ref{help1}) as $F_{j,r}^{w}(\alpha) = E[J_{c,j,w,r}^{UE}(\alpha)]$. Note that in a UE solution, it holds that $J_{c,j,w,r}^{UE}(\alpha) = J_{c,j,w^{\prime},r}^{UE}(\alpha), \forall w \neq w^{\prime}$.

At the UE, drivers make their own independent routing decisions. Therefore, in our formulation, given the assumption that truck drivers only know the probability distribution of the demand for the rest of the truck drivers and not the exact realization of it, their routing decisions  $\alpha_{j,w,r}^{UE}$ do not depend on the exact demand realization $c$. Given the aforementioned definitions, we can formulate the optimization problem through which we can calculate a UE solution as follows:
\begin{equation}
\label{usereq}
\begin{aligned}
& \underset{\alpha,\zeta}{\text{minimize}}
& & \lambda E[T_{tr}(\alpha)]+(1-\lambda)E[T^{mon}_{tr}(\alpha)] \\
& \text{subject to}
& & 0 \leq \alpha^{j}_{w,r} \perp F^w_{j,r}(\alpha)-\zeta^j_{w} \geq 0, \; \forall j,w,r \\
&&& \sum_{r \in R_j} \alpha^{j}_{w,r}=1, \; \forall j,w  
\end{aligned}
\end{equation}
where $\zeta_w^j$ is a set of free variables that are used in order to solve the equilibrium optimization problem (\ref{usereq}) and $F^w_{j,r}(\alpha)$ is given by (\ref{help1}). Additionally, the notation $\perp$ means that either $\alpha^{j}_{w,r}=0$ or $F^w_{j,r}(\alpha)-\zeta^j_{w}=0$ and finally, $\lambda$ is a weighting factor such that $\lambda\in[0,1]$. Therefore, in the equilibrium optimization problem (\ref{usereq}), among the possibly nonequivalent UE solutions, we are looking for the one that minimizes a weighted combination of the expected total travel time and the expected total monetary cost of the truck drivers. Setting $\lambda=1$ in the objective function of (\ref{usereq}), we can calculate the UE with the minimum expected total travel time.

Viewing the expected total travel time of the truck drivers as a uniformly weighted expected total cost and given the fact that $E[T^{mon}_{tr}(\alpha)]$ is equal to the expected total travel time of the truck drivers weighted by the corresponding VOT of each class $w$, the overall objective of (\ref{usereq}) can be expressed in cost units. The reasoning behind choosing the UE solution which minimizes the objective function of (\ref{usereq}) is the following. First, as also mentioned in \cite{RePEc:eee:transb:v:44:y::i:8-9:p:972-982}, $E[T_{tr}(\alpha)]$ has long been accepted as a standard index of system performance in a transportation context while $E[T^{mon}_{tr}(\alpha)]$ is a more appropriate system measure from an economic viewpoint. Second, we use the solution of (\ref{usereq}) as a benchmark for designing Pareto-improving pricing-and-routing schemes. Note that in order to create a Pareto-improving pricing scheme, i.e. a pricing scheme that can make everyone better-off compared to the UE, we first need to guarantee that the expected total travel time and the expected total monetary cost of the truck drivers using the proposed pricing scheme are lower than their best possible corresponding values at the UE.   

Recently, to study how close a real traffic scenario is to a UE, for the static traffic assignment problem \citep{patriksson}, \cite{Cabannes:2019:RRN:3350424.3325916} defined the average marginal regret as the expected time-saving drivers have in the network if they change their path to an optimal one. They proved that as the number of routing apps used is increased, the observed traffic assignment converges to a UE. The simulation results using real data for the whole Los Angeles network showed that the minimum travel time of a driver can be achieved whenever the ratio of GPS routing app users reaches 100\%. In this case, the network converges to the UE. In this work, we use the UE as a benchmark for our design since we may
expect that it is an optimistic version of the real world traffic
conditions where all truck drivers have updated information of
the traffic conditions through the use of GPS routing apps.

\subsection{System Optimum (SO).}
In a System Optimum (SO) solution, drivers are making routing decisions in a manner that contributes to the minimization of a social cost compared to the UE where they minimize their own individual travel time. Letting $E[T_{p}(\alpha)]$ denote the expected total travel time of the passenger vehicles in the network, we define the expected total travel time of the network as:
\begin{equation}
\label{total_trav_time}
\begin{aligned}
E[T_s(\alpha)] = E[T_{p}(\alpha)] + E[T_{tr}(\alpha)]
\end{aligned}
\end{equation}
Note that in a SO solution, the routing decisions of the truck drivers depend on the exact realization of the OD demands and therefore, we define the expected total monetary cost of the truck drivers as:
\begin{equation}
\label{with_c_mon_cost}
\begin{aligned}
E[T_{tr}^{mon}(\alpha)] = \sum_c \sum_{j=1}^{v} \sum_{w=1}^{N} \sum_{r \in R_j} p_c d_{j,w}^{c} \alpha_{w,r}^{c,j} s_w J_{c,j,w,r}(\alpha)
\end{aligned}
\end{equation}
where the main difference between (\ref{mon_cost}) and (\ref{with_c_mon_cost}) is the fact that in (\ref{with_c_mon_cost}), the routing decisions of the truck drivers $\alpha_{w,r}^{c,j}$ depend on the demand realization $c$. Note that $J_{c,j,w,r}$ has a similar definition as in (\ref{J_UE}), i.e.:
\begin{equation}
\label{J_mech}
    J_{c,j,w,r}(\alpha) = \sum_{l \in r} C_{lT}(X_{lp}, X_{lT}(\alpha))
\end{equation}
where $X_{lT}(\alpha)$ is given by:
\begin{equation}
\label{X_lT_again}
    X_{lT}(\alpha) = \sum_{j=1}^v \sum_{w=1}^N \sum_{r \in R_j: l \in r} d_{j,w}^{c} \alpha_{w,r}^{c,j}
\end{equation}
where the main difference between (\ref{X_lT_again}) and (\ref{proti}) is the fact that in (\ref{X_lT_again}), the routing decisions of the truck drivers $\alpha_{w,r}^{c,j}$ depend on the demand realization $c$. Since $J_{c,j,w,r}(\alpha)$ expresses travel time, it holds that $J_{c,j,w,r}(\alpha) = J_{c,j,w^{\prime},r}(\alpha), \forall w \neq w^{\prime}$. Using the aforementioned definitions, we calculate the SO solution of the network by solving the following optimization problem:  
\begin{equation}
\label{SO}
\begin{aligned}
& \underset{\alpha(\cdot)}{\text{minimize}}
& & \lambda (\mu E[T_{tr}(\alpha)] + (1-\mu)E[T_{p}(\alpha)])+(1-\lambda)E[T^{mon}_{tr}(\alpha)] \\
& \text{subject to}
& & \sum_{r \in R_j} \alpha^{c,j}_{w,r}=1, \; \forall c,j,w \\ 
&&& \alpha^{c,j}_{w,r}\geq 0, \; \forall c,j,w,r
\end{aligned}
\end{equation}
where $\lambda,\mu \in [0,1]$. In (\ref{SO}), we minimize a weighted combination of the expected total travel time of the truck drivers, their expected total monetary cost and the expected total travel time of the passenger vehicles. The reasoning behind the selection of this objective function is that even though we are providing routing suggestions only to the truck drivers, simultaneously, we want to make sure that the passenger vehicles will not be significantly affected. The parameter $\mu$ can be used to adjust the weight put on each category of vehicles. In all of the experiments, we set $\mu=0.9$.

\section{Pricing-and-Routing Schemes.}
In a UE solution, every driver makes his/her own individual routing decisions which leads to an inefficient road usage. On the other hand, in a SO solution, some drivers may benefit while some others may be harmed compared to the UE solution, providing no incentives to drivers to follow the SO solution in practice. In this section we study pricing-and-routing schemes that are Pareto-improving, i.e. they can make every user better-off compared to the UE while at the same time, they can drive the network as close as possible to the SO solution. Note that even though the proposed schemes are pricing-and-routing schemes, we often refer to them as pricing schemes.

We design a coordination mechanism that can be applied to truck drivers taking into account the user heterogeneity in their VOT. More specifically, the coordinator asks the truck drivers to declare their desired OD pair and additionally choose their VOT from a set of $N$ available options. We assume that a monitoring system is in place and thus the truck drivers may not lie about their declared OD pair. This assumption is not restrictive thanks to the advancements in GPS technologies. Additionally, we assume that GPS tracking is used to make sure that the truck drivers obey the routing suggestions. After collecting this information, the coordinator provides routing suggestions and additionally designs pricing schemes that are Pareto-improving and guarantee that every driver will have an incentive to truthfully declare his/her VOT while concurrently leading to a revenue-neutral (budget balanced) on average mechanism. This is in contrast with the previous literature studying pricing schemes, e.g. \citep{e430d1dc1f7a428d83c3c49d1a22bcfc, RePEc:eee:transe:v:54:y:2013:i:c:p:1-13}, that makes assumptions about the distribution that the user heterogeneity might follow. We should note that it is important to guarantee that a user will truthfully declare his/her VOT in order to avoid the exploitability of the designed mechanism. This is mainly because many truck drivers would be willing to declare a high VOT in order to be assigned to the fastest possible route. In the next two subsections, we design pricing-and-routing schemes that mathematically satisfy the property of truthfulness.

In the following sections, we propose two pricing-and-routing schemes, the Optimum Pricing Scheme (OPS) and the Approximately Optimum Pricing Scheme (AOPS). In both of these schemes, after asking the drivers to declare their OD pair $j$ and pick their VOT $w$ from a set of $N$ available options, we solve an optimization problem, we calculate a way to route the drivers into the network as well as the corresponding pricing scheme and finally, we randomly assign the truck drivers who declared the same OD pair $j$ and class $w$ into routes $r$.

\subsection{Optimum Pricing Scheme (OPS).}
Let $\pi^{c,j}_{w,r}$ be the payment (made or received) by a truck driver belonging to the class $w$ with an OD pair $j$ who follows route $r$ during demand realization $c$. We calculate the optimum way to route the truck drivers $\alpha^{*}$ as well as the the optimum pricing scheme $\pi^{*}$ by solving the following nonconvex optimization problem: 
\begin{equation}
\label{optim2}
\begin{aligned}
& \underset{\alpha(\cdot),\pi(\cdot)}{\text{minimize}}
& & \lambda (\mu E[T_{tr}(\alpha)] + (1-\mu)E[T_{p}(\alpha)])+(1-\lambda)E[T^{mon}_{tr}(\alpha)] \\
& \text{subject to}
& & \sum_c \sum_{r \in R_j} p_c  \alpha^{c,j}_{w,r} (J^{M,c,j}_{w,r} + \frac{1}{s_w} \pi^{c,j}_{w,r}) \leq \sum_c p_c A^{UE}_{c,j}, \; \forall j,w \\
&&& \sum_c \sum_{r \in R_j} p_c \alpha^{c,j}_{i,r} (J^{M,c,j}_{i,r} + \frac{1}{s_i} \pi^{c,j}_{i,r}) \leq \sum_c \sum_{r \in R_j} p_c \alpha^{c,j}_{k,r} (J^{M,c,j}_{k,r} + \frac{1}{s_i} \pi^{c,j}_{k,r}), \; \forall j,i,k \\
&&& \sum_{c} \sum_{j=1}^v \sum_{w=1}^N \sum_{r \in R_j} p_c d^{w}_{c,j} \alpha^{c,j}_{w,r} \pi^{c,j}_{w,r} = 0\\
&&& \sum_{r \in R_j} \alpha^{c,j}_{w,r}=1, \; \forall c,j,w \\ 
&&& \alpha^{c,j}_{w,r}\geq 0, \; \forall c,j,w,r
\end{aligned}
\end{equation}
where $J^{M,c,j}_{w,r}$ is the travel time of a truck driver belonging to class $w$ with a desired OD pair $j$ who follows route $r$ during the demand realization $c$ under the mechanism routing suggestions $M$ and is defined according to (\ref{J_mech}). Furthermore, $A_{c,j}^{UE}$ is the average travel time of a truck driver with OD pair $j$ during the demand realization $c$ at the UE and is given by the following equation:
\begin{equation}
\label{average_in_UE}
    A_{c,j}^{UE} = \sum_{r \in R_j} \alpha_{j,w,r}^{UE} J_{c,j,w,r}^{UE}
\end{equation}
Note that in (\ref{average_in_UE}), we omit the index $w$ from $A_{c,j}^{UE}$ since at the UE, the average travel times of truck drivers belonging to different classes of VOT are identical. Additionally, it holds that $\sum_c p_c A^{UE}_{c,j} = \sum_c p_c J^{UE}_{c,j,w,r}, \forall j,w,r$.
Based on the aforementioned definitions, the first constraint of (\ref{optim2}) guarantees that the expected cost (in time units) of the drivers at the time they make their decision is less than or equal to their corresponding cost at the UE (Pareto-improvement). Note that this definition of the Pareto-improvement property is similar to the one used in \cite{RePEc:eee:transb:v:44:y::i:8-9:p:972-982}. In our case, we extend this definition to the case where the demand is stochastic. The second constraint of (\ref{optim2}) guarantees that a truck driver which belongs to class $i$ and truthfully declares class $i$ to the coordinator will be better-off on average compared to the case where he/she originally belongs to class $i$ but declares class $k$ to the coordinator. Therefore, under the assumption that truck drivers are rational and will be constantly seeking to minimize their individual aggregated time (travel time + payments in time units), the second constraint of (\ref{optim2}) guarantees that every user will have an incentive to truthfully declare his/her VOT. Last, the third constraint of (\ref{optim2}) guarantees that the expected total payments made and received by the coordinator are equal to zero and therefore, the resulting mechanism satisfies the budget balanced on average property. Note that in the third constraint of (\ref{optim2}), we implicitly assume that the coordinator incurs no cost from operating the pricing-and-routing mechanism. This is a common assumption both in a game-theoretic context \citep{10.5555/1483085} and in a transportation context, e.g. \citep{RePEc:eee:transb:v:44:y::i:8-9:p:972-982}. At this point, note that the UE solution, where no pricing scheme is applied to the users, satisfies the constraints of (\ref{optim2}) and therefore, a solution to (\ref{optim2}) always exists. Hence, the desired properties hold regardless of the objective function that we choose to minimize in (\ref{usereq}).

At this point, let us comment on the Pareto-improvement and the truthfulness properties of OPS. As far as it concerns the Pareto-improvement property, as can be seen from the first constraint of (\ref{optim2}), on the left side of the inequality, we calculate the average over the routes $r$ and the expectation over the different demand realizations $c$. Theoretically, it is possible that for some demand realizations, some drivers are given a route with a higher total travel time (travel time + payments expressed in time units) compared to the average travel time at the UE. However, due to the fact that individual drivers only know the probability distribution of the demand and not the exact realization of it and hence they have incomplete information of the traffic conditions, they will be willing to participate in OPS since at the time they make their decision, their expected cost will be lower under the mechanism routing suggestions $M$ than in the UE. Additionally, if users repeatedly participate in OPS, then using randomization in order to assign the truck drivers who declared the same OD pair $j$ and class $w$ into routes $r$, the first constraint of (\ref{optim2}) guarantees that every truck driver will be better-off compared to the UE (Pareto-improvement). In Appendix~B, we prove the existence of pricing-and-routing schemes that guarantee an even stronger version of the Pareto-improvement property. More specifically, we prove that it is possible to design a scheme that guarantees that at every demand realization $c$ and for every route $r$ that a driver might be assigned to, his/her total travel time (travel time + payments expressed in time units) is going to be lower under the mechanism suggestions compared to the UE.

As far as it concerns the truthfulness property of OPS, as can be seen from the second constraint of (\ref{optim2}), we calculate the average over the routes $r$ and the expectation over the different demand realizations $c$. Again, using the fact that individual drivers only know the probability distribution of the demand and not the exact realization of it and hence they have incomplete information of the traffic conditions, it is guaranteed that they will be willing to truthfully declare their VOT since at the time they make their decision, their expected cost in the case they are truthful is lower than their corresponding cost in the case where they declared a different VOT than their actual one. In Appendix~B, we prove the existence of pricing-and-routing schemes that guarantee an even stronger version of the truthfulness property since we force the second constraint of (\ref{optim2}) to hold for every demand realization $c$.

Even though the schemes presented in Appendix B can guarantee even stronger versions of the Pareto-improvement and the truthfulness properties, we decided not to include them in our analysis for three main reasons. First, OPS is still sufficient to guarantee that the drivers will have an incentive to participate in the mechanism and truthfully declare their VOT since at the time they make their decision, their expected cost in the case they are truthful is lower than their corresponding cost in the case where they declared a different VOT than their actual one. Second, by forcing the second constraint of (\ref{optim2}) to hold for every possible demand realization $c$, the number of constraints would increase and the final scheme would become less computationally efficient. Finally, one can expect that the optimum solution would be less efficient since the size of the feasible region over which we optimize would be smaller.

In order to reduce the dimensionality of (\ref{optim2}), in the following subsection, we present an Approximately Optimum Pricing Scheme (AOPS) and show that we can assign routes to the drivers so that the proposed pricing scheme meets the desired goals.  

\subsection{Approximately Optimum Pricing Scheme (AOPS).}
For a given routing decision $\alpha$, let us define the following pricing scheme:
\begin{equation}
\label{payments}
\begin{aligned}
\pi^{AOPS}_{c,j,w,r} = \ &s_{w} (A^{UE}_{c,j}- J^{M,c,j}_{w,r}) + \frac{s_{w}}{\sum_{l=1}^N s_l} \frac{E\big[T^{mon,M}_{tr}\big] - E\big[T^{mon,UE}_{tr}\big]}{\sum_{j=1}^v d^{w}_{c,j}}
\end{aligned}
\end{equation}
The pricing scheme given by (\ref{payments}) initially makes each driver pay (or receive a payment) such that his/her travel time under the mechanism routing suggestions $J^{M,c,j}_{w,r}$ becomes equal to his/her average travel time at the UE $A^{UE}_{c,j}$. Then, after calculating the expected total monetary benefits of the truck drivers $E\big[T^{mon,M}_{tr}\big] - E\big[T^{mon,UE}_{tr}\big]$ obtained from the application of the mechanism, it distributes those benefits to the different classes proportionally to the VOT that each class has. Finally, each class benefits are uniformly shared among the truck drivers of the class.

Let us now formulate the following optimization problem:
\begin{equation}
\label{approx_optimization}
\begin{aligned}
& \underset{\alpha(\cdot)}{\text{minimize}}
& & \lambda (\mu E[T_{tr}(\alpha)] + (1-\mu)E[T_{p}(\alpha)])+(1-\lambda)E[T^{mon}_{tr}(\alpha)] \\
& \text{subject to}
& & E[T^{mon}_{tr}(\alpha)]\leq E\big[T^{mon,UE}_{tr}\big] \\
&&& H_{i,k}^{j}(\alpha) \leq N_{i,k}^{j}(\alpha), \; \forall j,i,k \\
&&& \sum_{r \in R_j} \alpha^{c,j}_{w,r}=1, \; \forall c,j,w \\
&&& \alpha^{c,j}_{w,r}\geq 0, \; \forall c,j,w,r
\end{aligned}
\end{equation}
where $E\big[T^{mon,UE}_{tr}\big]$ is the expected total monetary cost of the truck drivers at the UE and $H_{i,k}^{j}(\alpha)$ and $N_{i,k}^{j}(\alpha)$ are given by the following equations:
\begin{equation}
\label{cond2}
\begin{aligned}
H_{i,k}^{j}(\alpha) = \bigg(&1-\frac{s_k}{s_i}\bigg) \sum_c p_c A^{UE}_{c,j} + \frac{1}{\sum_{w=1}^N s_w} \sum_c p_c \frac{E[T^{mon,M}_{tr}] - E\big[T^{mon,UE}_{tr}\big]}{\sum_{j=1}^v d^i_{c,j}} 
\end{aligned}
\end{equation}
\begin{equation}
\label{cond3}
\begin{aligned}
N_{i,k}^{j}(\alpha) = \bigg(1-\frac{s_k}{s_i}\bigg) \sum_c \sum_{r \in R_j} p_c \alpha^{c,j}_{k,r} J^{M,c,j}_{k,r} + \frac{s_k}{s_i} \frac{1}{\sum_{w=1}^N s_w} \sum_c p_c \frac{E[T^{mon,M}_{tr}] - E\big[T^{mon,UE}_{tr}\big]}{\sum_{j=1}^v d^k_{c,j}}
\end{aligned}
\end{equation}
Note that a solution to the optimization problem (\ref{approx_optimization}) always exists since the UE satisfies all of its constraints. Let us call the optimum solution of the optimization problem described by (\ref{approx_optimization})-(\ref{cond3}) as $\alpha_{AOPS}^{*}$. Now, we are ready to state the following theorem.

\begin{theorem}
The pair ($\alpha_{AOPS}^{*}$, $\pi^{AOPS}_{c,j,w,r}$) makes the expected cost of the drivers at the time they make their decision to be less than or equal to their corresponding cost at the UE, guarantees that every user will have an incentive to truthfully declare his/her VOT and leads to a budget balanced on average mechanism.
\end{theorem}
{{ \bf Proof.} To prove the statement of the theorem, we equivalently prove that $\pi^{AOPS}_{c,j,w,r}$ is Pareto-improving, guarantees that every user will have an incentive to truthfully declare his/her VOT and creates a budget balanced on average mechanism if and only if the first and the second constraint of (\ref{approx_optimization}) together with (\ref{cond2})-(\ref{cond3}) hold. Note that a user will be better-off compared to the UE if the first constraint of (\ref{optim2}) holds. Therefore, substituting (\ref{payments}) into the first constraint of (\ref{optim2}), we get:
\begin{equation*}
\begin{aligned}
\sum_c \sum_{r \in R_j} p_c \alpha^{c,j}_{w,r} \bigg(J^{M,c,j}_{w,r} + A^{UE}_{c,j} - J^{M,c,j}_{w,r} + \frac{1}{\sum_{l=1}^N s_l} \frac{E[T^{mon,M}_{tr}] - E\big[T^{mon,UE}_{tr}\big]}{\sum_{j=1}^v d^{w}_{c,j}} \bigg) \leq \sum_c p_c A^{UE}_{c,j} \Leftrightarrow \\
 \Leftrightarrow \frac{1}{\sum_{l=1}^N s_l} \big(E[T^{mon,M}_{tr}] - E\big[T^{mon,UE}_{tr}\big]\big) \sum_c \frac{p_c}{\sum_{j=1}^v d^w_{c,j}} \leq 0
\end{aligned}
\end{equation*}
which holds true if and only if $E[T^{mon}_{tr}(\alpha)]\leq E\big[T^{mon,UE}_{tr}\big]$ which is equivalent to the first constraint of (\ref{approx_optimization}). Additionally, a user will have an incentive to truthfully declare his/her VOT if the second constraint of (\ref{optim2}) holds. Therefore, substituting (\ref{payments}) into the second constraint of (\ref{optim2}), we get:
\begin{equation*}
\begin{aligned}
\sum_c \sum_{r \in R_j} p_c \alpha^{c,j}_{i,r} \bigg(J^{M,c,j}_{i,r} + A^{UE}_{c,j} - J^{M,c,j}_{i,r} + \frac{1}{\sum_{l=1}^N s_l} \frac{E[T^{mon,M}_{tr}] - E\big[T^{mon,UE}_{tr}\big]}{\sum_{j=1}^v d^{i}_{c,j}} \bigg) \leq \\ \leq
\sum_c \sum_{r \in R_j} p_c \alpha^{c,j}_{k,r} \bigg(J^{M,c,j}_{k,r} + \frac{s_k}{s_i} \bigg( A^{UE}_{c,j} - J^{M,c,j}_{k,r} + \frac{1}{\sum_{l=1}^N s_l} \frac{E[T^{mon,M}_{tr}] - E\big[T^{mon,UE}_{tr}\big]}{\sum_{j=1}^v d^{k}_{c,j}} \bigg)\bigg) \Leftrightarrow \\
\Leftrightarrow  
\bigg(1-\frac{s_k}{s_i}\bigg) \sum_c p_c A^{UE}_{c,j} + \frac{1}{\sum_{l=1}^N s_l} \sum_c p_c \frac{E[T^{mon,M}_{tr}] - E\big[T^{mon,UE}_{tr}\big]}{\sum_{j=1}^v d^i_{c,j}} \leq \\
\leq  \bigg(1-\frac{s_k}{s_i}\bigg) \sum_c \sum_{r \in R_j} p_c \alpha^{c,j}_{k,r} J^{M,c,j}_{k,r} + \frac{s_k}{s_i} \frac{1}{\sum_{l=1}^N s_l} \sum_c p_c \frac{E[T^{mon,M}_{tr}] - E\big[T^{mon,UE}_{tr}\big]}{\sum_{j=1}^v d^k_{c,j}} 
\end{aligned}
\end{equation*}
where the last inequality is equivalent to the second constraint of (\ref{approx_optimization}). Last, a mechanism is budget balanced on average if the third constraint of (\ref{optim2}) holds. Substituting (\ref{payments}) into the third constraint of (\ref{optim2}), we get:
\begin{equation*}
\begin{aligned}
\sum_{c} \sum_{j=1}^v \sum_{w=1}^N \sum_{r \in R_j} p_c d^{w}_{c,j} \alpha^{c,j}_{w,r} \bigg(s_{w} (A^{UE}_{c,j}- J^{M,c,j}_{w,r}) + \frac{s_{w}}{\sum_{l=1}^N s_l} \frac{E[T^{mon,M}_{tr}] - E\big[T^{mon,UE}_{tr}\big]}{\sum_{j=1}^v d^{w}_{c,j}} \bigg) \\
= E\big[T^{mon,UE}_{tr}\big] - E\big[T^{mon,M}_{tr}\big] + \sum_c \sum_{j=1}^v \sum_{w=1}^N p_c d^{w}_{c,j}
\frac{s_{w}}{\sum_{l=1}^N s_l} \frac{E[T^{mon,M}_{tr}] - E\big[T^{mon,UE}_{tr}\big]}{\sum_{j=1}^v d^{w}_{c,j}} \sum_{r \in R_j} \alpha^{c,j}_{w,r}
\\
= E\big[T^{mon,UE}_{tr}\big] - E\big[T^{mon,M}_{tr}\big] + \big(E\big[T^{mon,M}_{tr}\big] - E\big[T^{mon,UE}_{tr}\big] \big) \bigg( \sum_c p_c \frac{1}{\sum_{l=1}^{N} s_l} \sum_{w=1}^{N} s_w \frac{1}{\sum_{j=1}^v d^{w}_{c,j}} \sum_{j=1}^{v} d_{c,j}^{w} \bigg)
\\
= E\big[T^{mon,UE}_{tr}\big] - E\big[T^{mon,M}_{tr}\big] + E\big[T^{mon,M}_{tr}\big] - E\big[T^{mon,UE}_{tr}\big] = 0
\end{aligned}
\end{equation*}
where in the second equality, we used the feasibility constraint which is given by the third constraint of (\ref{approx_optimization}). Since the UE satisfies the constraints of (\ref{approx_optimization}), a solution to (\ref{approx_optimization}) always exists. Therefore, we have proved that by solving the optimization problem described by (\ref{approx_optimization})-(\ref{cond3}), we can calculate $\alpha_{AOPS}^{*}$ such that the pricing scheme $\pi^{AOPS}_{c,j,w,r}$ given by (\ref{payments}) satisfies all the statements of Theorem~1 and this concludes the proof.\hfill \qed}

Theorem 1 states that one can get a sub-optimal solution to the original optimization problem (\ref{optim2}) by solving the optimization problem described by (\ref{approx_optimization})-(\ref{cond3}) in order to assign routes to the truck drivers and by subsequently applying the pricing scheme given by (\ref{payments}). We call this method Approximately Optimum Pricing Scheme (AOPS). The main advantage of this approach is the fact that we significantly reduce the dimensionality of the problem by calculating a pricing scheme using a simple algebraic equation. As we will later show experimentally, AOPS achieves a significant improvement compared to the UE and provides a solution close to the SO.  

In the first part of the proof of Theorem 1, we proved that the condition $E[T^{mon}_{tr}(\alpha)]\leq E\big[T^{mon,UE}_{tr}\big]$ is a necessary and sufficient condition to guarantee that the pricing scheme (\ref{payments}) is Pareto-improving. For the deterministic demand case, a similar property was also proved in \cite{RePEc:eee:transb:v:44:y::i:8-9:p:972-982} for the case of revenue refunding schemes.

At this point, let us comment on the Pareto-improvement and the truthfulness properties of AOPS. Similar to OPS, it is theoretically possible that for some demand realizations, some drivers are given a route with a higher total travel time (travel time + payments expressed in time units) compared to the average travel time at the UE. However, due to the fact that individual drivers only know the probability distribution of the demand and not the exact realization of it and hence they have incomplete information of the traffic conditions, they will be willing to participate in AOPS since at the time they make their decision, their expected cost will be lower under the mechanism routing suggestions $M$ than in the UE. Additionally, if users repeatedly participate in AOPS, then using randomization in order to assign the truck drivers who declared the same OD pair $j$ and class $w$ into routes $r$, it is guaranteed that every truck driver will be better-off compared to the UE (Pareto-improvement).
As far as it concerns the truthfulness property, AOPS guarantees that the users will be willing to truthfully declare their VOT since at the time they make their decision, their expected cost in the case they are truthful is lower than their corresponding cost in the case where they declared a different VOT than their actual one. In Appendix~B, we prove the existence of pricing-and-routing schemes that can guarantee even stronger versions of the Pareto-improvement and the truthfulness properties. However, we decided not to include them in our analysis for the reasons explained in Section 3.1.

In the following subsection, we present a class-anonymous pricing scheme, namely, the Congestion Pricing with Uniform Revenue Refunding (CPURR) scheme. This is going to be used later as one of the baselines in the simulation experiments.

\subsection{Congestion Pricing with Uniform Revenue Refunding (CPURR).}
Under a congestion pricing scheme, each driver is assigned a fee depending on the OD pair and the route he/she follows. \cite{RePEc:eee:transb:v:44:y::i:8-9:p:972-982} proposed to combine Congestion Pricing with a Uniform Revenue Refunding (CPURR) scheme, i.e. the fees collected from congestion pricing are uniformly distributed among the participant drivers irrespective of the class they belong. Therefore, the whole scheme is class-anonymous. 

The original CPURR scheme proposed in \cite{RePEc:eee:transb:v:44:y::i:8-9:p:972-982} is OD-based. However, since our approach focuses on route-based pricing schemes, in this paper, we calculate the route-based equivalent of the CPURR scheme. Therefore, the CPURR scheme can be calculated by solving the following optimization problem with complementarity constraints: 
\begin{equation}
\label{cpurr}
\begin{aligned}
& \underset{\alpha,\pi,\zeta}{\text{minimize}}
& & \lambda (\mu E[T_{tr}(\alpha)] + (1-\mu)E[T_{p}(\alpha)])+(1-\lambda)E[T^{mon}_{tr}(\alpha)] \\
& \text{subject to}
& & 0 \leq \alpha^{j}_{w,r} \perp F^{CP}_{j,w,r}(\alpha,\pi)-\zeta^j_{w} \geq 0, \; \forall j,w,r\\
&&& \sum_{r \in R_j} \alpha^{j}_{w,r}=1, \; \forall j,w \\ 
&&& \sum_{c} \sum_{j=1}^v \sum_{w=1}^N \sum_{r \in R_j} p_c d^{w}_{c,j} \alpha^{j}_{w,r} \pi^{j}_{r} = 0
\end{aligned}
\end{equation}
where $\zeta_w^j$ is a set of free variables that are used in order to solve the equilibrium optimization problem (\ref{cpurr}) and $F^{CP}_{j,w,r}(\alpha,\pi)$ is given by the following equation:
\begin{equation}
F^{CP}_{j,w,r}(\alpha,\pi)=\sum_c p_c \bigg(J^{CP}_{c,j,w,r}(\alpha)+\frac{1}{s_w}\pi_r^j\bigg) 
\end{equation}
where $J^{CP}_{c,j,w,r}(\alpha)$ is defined according to (\ref{J_UE}). Note that under a congestion pricing scheme, the network users make their own individual routing decisions while taking into account the fees corresponding to each route. Since in our model the truck drivers only know the probability distribution of the demand for the rest of the truck drivers and not its exact realization, the way that the drivers choose their routes $\alpha^{j}_{w,r}$ does not depend on the exact realization $c$. Additionally, note that the variables $\pi^{j}_{r}$ corresponding to the pricing scheme do not depend on $w$ and $c$. The independence of $w$ can be justified by the fact that the CPURR scheme is class-anonymous. On the other hand, the coordinator of the CPURR scheme who is responsible for assigning fees to each route and then uniformly distribute the collected revenue back to the participant drivers, could design a pricing scheme that depends on the exact realization $c$ of the demand. However, even in that case, since none of the constraints of (\ref{cpurr}) depends on $c$, the optimum solution of (\ref{cpurr}) would not change.

In the next section, CPURR is used as one of the baselines in the simulation experiments. There are several reasons that justify this choice. First, as mentioned earlier, congestion pricing is one of the most widely studied methods for addressing the problem of the inefficiency between the UE and the SO solutions. Additionally, CPURR is a class-anonymous pricing scheme. This fact allows us to directly demonstrate the benefits of applying class-specific (VOT-based) pricing schemes such as OPS and AOPS.

\section{Experimental Results.}
In this section, we demonstrate our approach conducting simulation experiments based on the Sioux Falls network \citep{siouxfalls}. The Sioux Falls network consists of 24 nodes and 76 links and constitutes a benchmark in the transportation research field. The experimental results section is divided into four subsections: In the first subsection, we test the sensitivity of the solutions of the proposed methods to the initial conditions and to the number of routes considered per OD pair. Using these results, in the second subsection, we compare the UE, SO, OPS, AOPS and the CPURR scheme in terms of the expected total travel time of the truck drivers $E[T_{tr}]$, their expected total monetary cost $E[T_{tr}^{mon}]$, the expected total travel time of the network $E[T_{s}]$ and the total objective value by varying the weighting factor $\lambda$ and we show that the VOT-based pricing schemes (OPS and AOPS) outperform class-anonymous pricing schemes like CPURR. Subsequently, in the third subsection, we experimentally show that OPS and AOPS can be efficiently used under both a deterministic and a stochastic demand scenario since they outperform the CPURR scheme in terms of the expected total monetary cost of the truck drivers $E[T_{tr}^{mon}]$ and the total objective value by concurrently achieving a superior solution compared to the UE in both scenarios. In the last subsection, we increase the number of OD pairs that the truck drivers use in the Sioux Falls network and we experimentally show that both OPS and AOPS need a lower computational time compared to the CPURR scheme. Additionally, we show that AOPS remains computationally tractable even in the case where a large number of OD pairs is used by the truck drivers. For all of the experiments, the fmincon optimization solver implemented in the MATLAB Optimization Toolbox \citep{matlab} was used. Since fmincon solves optimization problems with local optimality guarantees, in this section, we compare local minima between the approaches.

\subsection{Robustness in the Sioux Falls Network}
 In our experiments, we assumed that the cost of each route corresponds to travel time and can be described by a Bureau of Public Roads (BPR) function \citep{sheffi} of the form:
\begin{equation*}
    \begin{aligned}
    C_{lT}(X_{lp},X_{lT}) = \epsilon_a + \epsilon_b \Bigg(\frac{X_{lp}+3X_{lT}}{\epsilon_c}\Bigg)^4
    \end{aligned}
\end{equation*}
where $\epsilon_a, \epsilon_b$ and $\epsilon_c$ are constants and their values were chosen similar to the ones adopted in \cite{papadopoulos2019incentives}\footnote{These values can be found in this \href{https://www.dropbox.com/s/rjqmac6ddqi9vd5/Parameters_of_the_BPR_function_Part_C.txt?dl=0}{link}.}. We further assumed that the number of passenger vehicles at each link of the Sioux Falls network is constant\footnote{These values can be found in this \href{https://www.dropbox.com/s/j6xfgm6wdfe9p26/Passenger_Vehicles.txt?dl=0}{link}.}. These numbers have been calculated by solving an equilibrium assignment problem for the passenger vehicles. To retain computational tractability, we further assumed that there are 6 available OD pairs for the truck drivers, namely $(O_1,D_7), (O_1,D_{11}), (O_{10},D_{11}), (O_{10},D_{20}), (O_{15},D_{5})$ and $(O_{24},D_{10})$ and that the truck drivers choose their class among two available options with VOT $s_1 = 200 \frac{\$}{hr}$ and $s_2 = 50 \frac{\$}{hr}$, respectively. The demand of the truck drivers was assumed to be stochastic and take one of the following two equiprobable values:
\begin{align*}
d_1 =
\begin{bmatrix}
3 & 4.5 & 6 & 3 & 14 & 3.6 \\
1 & 2.8 & 5.4 & 7 & 9 & 2
\end{bmatrix},
d_2 =
\begin{bmatrix}
5 & 1.8 & 3.9 & 15 & 6.4 & 2.4 \\
6 & 5.5 & 1.8 & 6.5 & 11 & 6
\end{bmatrix}
\end{align*}
where each column of $d_1$ and $d_2$ corresponds to the demand of truck drivers for each OD pair and each row denotes a different class of users. The ratio of trucks in the network was $5.16\%$ in the case of $d_1$ and $5.96\%$ in the case of $d_2$. The values of the weighting factors were chosen to be $\lambda = \mu = 0.9$.

\subsubsection{Sensitivity to Initial Conditions.}
To test the sensitivity of the solutions of the optimization problems of the proposed methods to different initial conditions, in this section, we designed the following test scenario. We randomly initialized the initial conditions for the UE problem. After solving the UE problem, we used its solution as an initial condition to the rest of the methods, i.e. SO, OPS, AOPS and CPURR.

In our experiments, we assumed that the truck drivers follow the 10 least congested routes per OD pair. To calculate the least congested routes for each OD pair, we followed the following procedure. Before assigning any truck driver into the network ($X_{lT}=0$), we calculated the cost of each route by substituting the number of passenger vehicles at each road segment $l$ into the BPR function. Therefore, for each OD pair, we could calculate the least congested routes, i.e. the routes with the lowest travel time as defined by the BPR function when $X_{lT}=0$.

For the UE problem described by (\ref{usereq}), we picked the initial values for the set of free variables $\zeta_{w}^{j}$ according to a uniform distribution $\cal U[0,100]$. The values for the variables $\alpha_{w,r}^{j}$ were picked according to a uniform distribution $\cal U[0,1]$. The interior point method provided in the MATLAB Optimization Toolbox \citep{matlab} was used in order to solve (\ref{usereq}). Subsequently, this solution was used as an initial condition for the optimization problems of the other methods. Note that this solution is a feasible solution for all of the methods. This experiment was repeated 10 times. In Table~\ref{solution_robustness}, we report the minimum, the maximum, the mean, the standard deviation and the Coefficient of Variation (CV) values for different metrics for the UE, SO, OPS, AOPS and the CPURR methods.

The results of Table~\ref{solution_robustness} demonstrate that both OPS and AOPS approaches outperform the CPURR scheme and can approach the SO solution. Additionally, the solutions of OPS and AOPS are shown to be robust since the standard deviation values of these methods are significantly lower compared to the ones of the CPURR scheme.

\begin{table*}[h]
\begin{center}
\begin{adjustbox}{max width=\textwidth}
\begin{tabular}{cc|c|c|c|c|c|}
\\
\cline{3-7} 
Method&Metric&\textbf{Min}&\textbf{Max}&\textbf{Mean}&\textbf{Std}&\textbf{CV (\%)}\\
\hline
\multirow{5}{*}{{{UE}}}&$E[T_{tr}]$&20043.3&20190.7&20117.8&48.6&0.242\\
&$E[T_{tr}^{mon}]$&43235.8&43567.8&43414.8&128.0&0.295\\
&$E[T_{s}]$&330235.8&330826.2&330534.7&192.7&0.058\\
&$E[T_{p}]$&310192.5&310635.5&310416.9&144.1&0.046\\
&Objective&48476.0&48667.1&48574.4&64.5&0.133\\
\hline
\multirow{5}{*}{{{SO}}}&$E[T_{tr}]$&19682.2&19683.1&19682.5&0.3&0.002\\
&$E[T_{tr}^{mon}]$&41727.8&41777.7&41741.3&13.3&0.032\\
&$E[T_{s}]$&329109.4&329129.5&329125.6&5.7&0.002\\
&$E[T_{p}]$&309426.6&309447.3&309443.1&5.8&0.002\\
&Objective&47966.2&47969.2&47966.9&0.8&0.002\\
\hline
\multirow{5}{*}{{{OPS}}}&$E[T_{tr}]$&19680.8&19683.3&19682.9&0.7&0.004\\
&$E[T_{tr}^{mon}]$&41729.1&41791.0&41759.4&20.9&0.050\\
&$E[T_{s}]$&329111.6&329130.3&329121.2&6.9&0.002\\
&$E[T_{p}]$&309428.5&309447.1&309438.3&7.0&0.002\\
&Objective&47966.6&47970.5&47968.5&1.3&0.003\\
\hline
\multirow{5}{*}{{{AOPS}}}&$E[T_{tr}]$&19698.8&19729.4&19716.8&10.1&0.051\\
&$E[T_{tr}^{mon}]$&41737.2&41821.7&41783.2&28.5&0.068\\
&$E[T_{s}]$&329274.0&329423.7&329378.9&41.9&0.013\\
&$E[T_{p}]$&309575.1&309719.1&309662.0&36.6&0.012\\
&Objective&47997.3&48027.0&48018.6&8.3&0.017\\
\hline
\multirow{5}{*}{{{CPURR}}}&$E[T_{tr}]$&19853.4&20070.8&19889.2&62.2&0.313\\
&$E[T_{tr}^{mon}]$&42632.8&43195.7&42737.9&162.0&0.379\\
&$E[T_{s}]$&329554.6&330565.8&329777.7&270.8&0.082\\
&$E[T_{p}]$&309681.4&310494.9&309888.5&210.2&0.068\\
&Objective&48227.5&48521.5&48274.0&84.9&0.176\\
\hline
\end{tabular}
\end{adjustbox}
\end{center}
\caption{\label{solution_robustness}The minimum, the maximum, the mean, the standard deviation and the coefficient of variation (CV) values for different metrics for the UE, SO, OPS, AOPS and CPURR methods after 10 simulation runs with different initial conditions. The number of OD pairs was 6 and we considered the 10 least congested routes per OD pair.}
\end{table*}

To further explore the robustness of the solutions to the optimization problems of different methods, we conducted some additional experiments. More specifically, we considered solving all the optimization problems with random initial conditions instead of using the UE solution as an initial condition. Note that since the interior point method requires an initial condition that is close to primal and dual feasibility, the random initial conditions for this method were created by adding a random component to the UE solution. Furthermore, we explored solving the problems with the Sequential Quadratic Programming (SQP) method provided in the MATLAB Optimization Toolbox. The experimental results showed that both SO, OPS and AOPS can be efficiently solved under various initial conditions both with the interior point and the SQP methods. On the other hand, the CPURR problem could not achieve a significantly better solution than the UE in the case where the SQP optimization method was used. Additionally, in case the CPURR method was initialized with a random initial condition, the interior point method would fail to converge.

Last, we need to mention that for all of the experiments of this section and for the rest of the paper, the step tolerance used in the fmincon solver was $10^{-6}$. For UE, SO, OPS and AOPS, the constraint tolerance was set to $10^{-6}$. However, due to the difficulties on calculating the optimum solution of the CPURR problem, the constraint tolerance for this problem was set to $5 \cdot 10^{-5}$.

For the rest of this paper, the solutions of the optimization problems have been calculated using the following procedure.
We randomly initialized the initial conditions for the UE problem. After solving the UE problem, we used its solution as an initial condition to the rest of the methods, i.e. SO, OPS, AOPS and CPURR. In the next subsection, we test the the sensitivity of the solution to the number of routes considered per OD pair in the Sioux Falls network.

\subsubsection{Sensitivity to the Number of Routes.}
In this section, we test the sensitivity of the solutions of the optimization problems of the proposed methods to different number of routes considered per OD pair. In Table~\ref{Number_of_routes}, we calculate the expected total truck travel time, the expected total truck monetary cost, the expected total network time and the expected total travel time of the passenger vehicles of the UE, SO, OPS, AOPS and CPURR for the case where the truck drivers use 6 OD pairs of the Sioux Falls network, considering the 5, 10 and 15 least congested routes per OD pair. The rest of the parameters were chosen similar to the ones used in Section 4.1. 

\begin{table*}[h]
\begin{center}
\begin{adjustbox}{max width=\textwidth}
\begin{tabular}{cc|c|c|c|c|c|}
\\
\cline{3-7} 
\# of routes&Metric&\textbf{UE}&\textbf{SO}&\textbf{OPS}&\textbf{AOPS}&\textbf{CPURR}\\
\hline
\multirow{4}{*}{{{5}}}&$E[T_{tr}]$&21680.4&21582.6&21581.6&21580.4&21598.7\\
&$E[T_{tr}^{mon}]$&46870.0&46250.6&46254.7&46307.6&46667.7\\
&$E[T_{s}]$&334508.2&334258.0&334262.3&334256.2&334279.4\\
&$E[T_{p}]$&312827.8&312675.4&312680.7&312675.8&312680.7\\
\hline
\multirow{4}{*}{{{10}}}&$E[T_{tr}]$&20043.3&19682.2&19683.1&19698.8&19857.9\\
&$E[T_{tr}^{mon}]$&43235.8&41740.1&41745.0&41794.7&42633.1\\
&$E[T_{s}]$&330235.8&329127.6&329124.7&329274.0&329685.7\\
&$E[T_{p}]$&310192.5&309445.3&309441.0&309575.1&309827.7\\
\hline
\multirow{4}{*}{{{15}}}&$E[T_{tr}]$&19521.8&18982.9&18982.8&19097.5&19400.1\\
&$E[T_{tr}^{mon}]$&41996.2&40168.8&40177.4&40439.1&41683.9\\
&$E[T_{s}]$&326222.4&325966.3&325970.7&325971.5&325925.8\\
&$E[T_{p}]$&306700.6&306983.4&306987.9&306828.3&306525.8\\
\hline
\end{tabular}
\end{adjustbox}
\end{center}
\caption{\label{Number_of_routes}The expected total truck travel time $E[T_{tr}]$, the expected total truck monetary cost $E[T_{tr}^{mon}]$, the expected total network time $E[T_{s}]$ and the expected total travel time of the passenger vehicles $E[T_p]$ of the UE, SO, OPS, AOPS and CPURR in the case where the truck drivers follow 6 OD pairs, considering the 5, 10 and 15 least congested routes per OD pair.}
\end{table*}

As can be observed from the results presented in Table~\ref{Number_of_routes}, the more routes we consider per OD pair, the lower the values we can achieve in all four evaluation metrics. However, as the number of routes considered per OD pair increases, the computational time also increases. Our simulations showed that by considering 10 routes per OD pair, we can achieve a good balance between network efficiency and computational time in the Sioux Falls network. Therefore, for the rest of the paper, we only consider the 10 least congested routes per OD pair.

\subsection{The Effect of the Weighting Factor $\lambda$.}
In this section, we conduct additional simulation experiments in order to demonstrate the effect of the weighting factor $\lambda$ in the solutions of the UE, SO, OPS, AOPS and CPURR approaches and to show that the VOT-based pricing schemes (OPS and AOPS) outperform both the CPURR (a class-anonymous pricing scheme) and the UE solution. For the Sioux Falls network, the configurations were chosen similar to the ones used in Section 4.1. 

\begin{figure}[h]
\begin{center}
\subfloat[Expected total truck travel time.]{\includegraphics[width=.45\textwidth]{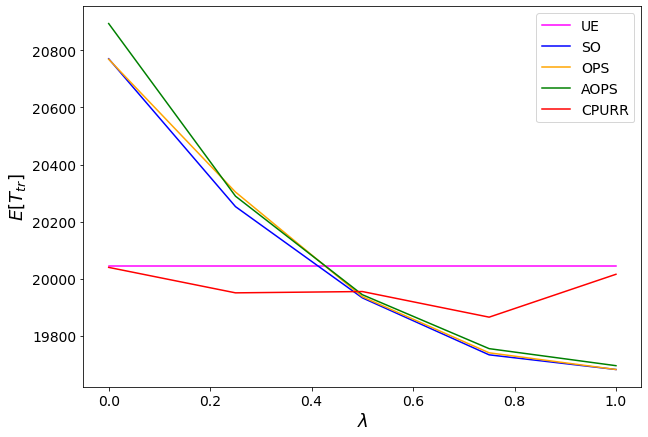}}
\subfloat[Expected total truck monetary cost.]{\includegraphics[width=.45\textwidth]{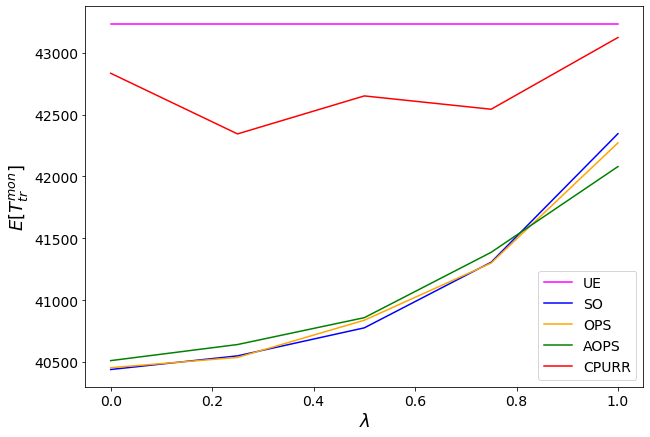}}\\
\subfloat[Expected total network time.]{\includegraphics[width=.45\textwidth]{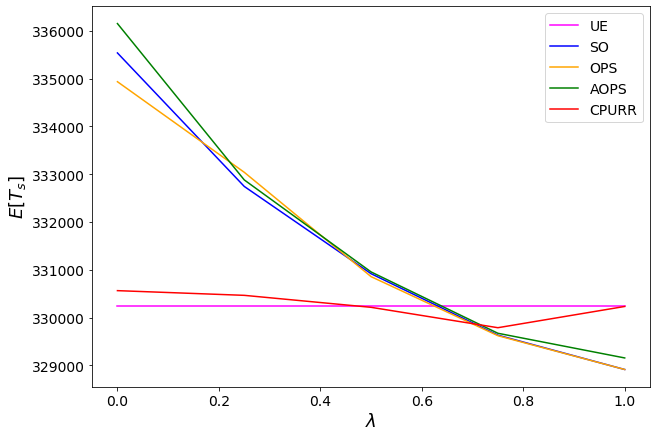}}
\subfloat[Total objective value.]{\includegraphics[width=.45\textwidth]{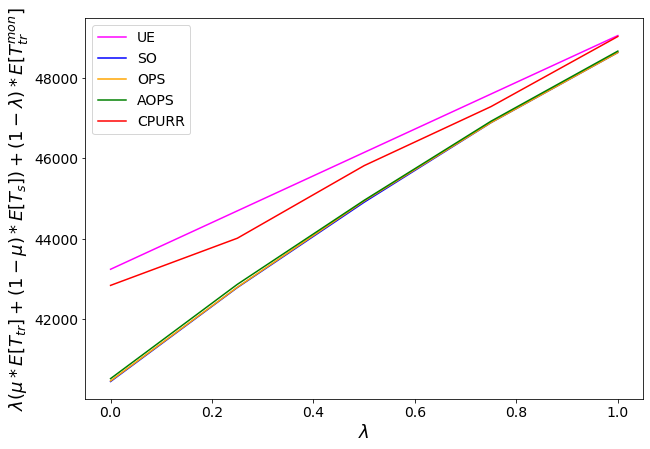}}
\caption{(a) The expected total truck travel time, (b) the expected total truck monetary cost, (c) the expected total network time and (d) the total objective value of the UE (magenta), SO (blue), OPS (orange), AOPS (green) and CPURR (red) for different values of the weighting factor $\lambda$.}
\label{multiplots}
\end{center}
\end{figure}

In Figure~\ref{multiplots}(a-d), we plot the expected total travel time of the truck drivers $E[T_{tr}]$, their expected total monetary cost $E[T_{tr}^{mon}]$, the expected total travel time of the network (passenger vehicles + trucks) and the total objective value, respectively, for different values of the weighting factor $\lambda \in [0,1]$.

As can be observed in Figure~\ref{multiplots}(a), as the value of the weighting factor $\lambda$ increases, the expected travel time of the truck drivers $E[T_{tr}]$ decreases for both SO, OPS and AOPS solutions. It is worth mentioning that for all values of $\lambda$, the OPS solution closely follows the SO solution. Additionally, the AOPS solution can significantly decrease $E[T_{tr}]$ compared to the CPURR solution, especially for $\lambda > 0.5$. 

In Figure~\ref{multiplots}(b), we observe that the expected total monetary cost of the truck drivers $E[T_{tr}^{mon}]$ increases as the value of the weighting factor $\lambda$ increases for both SO, OPS and AOPS solutions. On the other hand, we observe that $E[T_{tr}^{mon}]$ does not significantly change for the CPURR solution. Note also that AOPS has a smaller increase rate compared to the SO and the OPS solutions. This can be explained by the fact that as $\lambda$ increases, SO and OPS put more emphasis on minimizing the expected total travel time of the network while on the other hand, AOPS applies a pricing scheme during which the expected total monetary benefits $E[T_{tr}^{mon,M}]- E[T_{tr}^{mon,UE}]$, are shared to the users proportionally to the VOT of the class they belong. Therefore, a truck driver with higher VOT will get reimbursed with a bigger amount of money compared to a truck driver with lower VOT in case both of the drivers are assigned to a slower route. This behavior makes AOPS better contribute to the minimization of the expected total monetary cost of the truck drivers.    

In Figure~\ref{multiplots}(c), it is shown that as the value of the weighting factor $\lambda$ increases, the expected total travel time of the network $E[T_{s}]$ decreases. Furthermore, it can be observed that OPS makes $E[T_{s}]$ closely follow its corresponding value at the SO solution while at the same time, both OPS and AOPS can reduce the expected total travel time of the network compared to the CPURR scheme for $\lambda > 0.75$.

Finally, in Figure~\ref{multiplots}(d), we plot the total objective value $\lambda (\mu E[T_{tr}] + (1-\mu) E[T_{p}]) + (1-\lambda)E[T_{tr}^{mon}]$ and we observe that OPS closely approaches the SO solution. Additionally, AOPS provides a solution close to the SO and the OPS solutions, while constantly outperforming the CPURR scheme.

\subsection{Deterministic vs Stochastic Demand Scenario.}
In this subsection, we show that OPS and AOPS can be efficiently used under both a deterministic and a stochastic demand scenario. The network configurations were chosen similar to the ones used in Section 4.1.

For the purpose of this simulation experiment, we ran two distinct scenarios. In the first scenario (deterministic), we assumed that the truck drivers know the exact number of both the passenger vehicles and the rest of the trucks in the network. More specifically, the demand vector for the truck drivers was assumed to be:
\begin{equation*}
d = \begin{bmatrix}
3 & 4.5 & 6 & 3 & 14 & 3.6 \\
1 & 2.8 & 5.4 & 7 & 9 & 2
\end{bmatrix}
\end{equation*}
The ratio of truck drivers in the network was $5.16\%$ for this scenario. In the second scenario (stochastic), we assumed that the truck drivers know the number of passenger vehicles in the network and only the probability distribution of the demand for the rest of the truck drivers. More specifically, the demand of the truck drivers was assumed to take one of the two equiprobable values, $d_1$ (ratio of trucks: $5.16\%$) and $d_2$ (ratio of trucks: $5.96\%$), as given in Section 4.1. The values of the weighting factors were chosen to be $\lambda = \mu = 0.9$. 

\begin{figure}[h]
\begin{center}
\subfloat[Expected total truck monetary cost.]{\includegraphics[width=.45\textwidth]{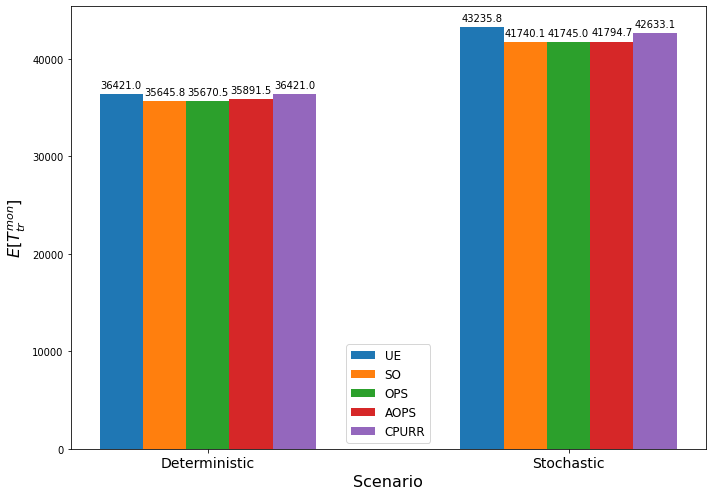}}
\subfloat[Objective value.]{\includegraphics[width=.45\textwidth]{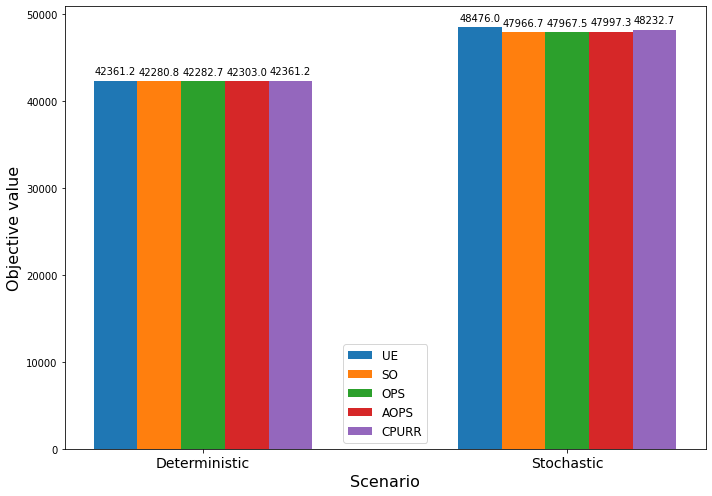}}
\caption{(a) The expected total truck monetary cost and (b) the total objective value of the UE (blue), SO (orange), OPS (green), AOPS (red) and CPURR (purple) for different demand scenarios.}
\label{deterministic_vs_stochastic}
\end{center}
\end{figure}

In Figure~\ref{deterministic_vs_stochastic}, we show the expected total truck monetary cost and the total objective value for the two scenarios. As can be observed from the results of Figure~\ref{deterministic_vs_stochastic}, OPS and AOPS approach the SO solution and outperform the UE and the CPURR scheme both in the deterministic and in the stochastic demand scenarios. Additionally, the benefits of applying OPS and AOPS rather than the CPURR scheme become greater in the stochastic demand scenario. The latter result was expected since OPS and AOPS provide routing suggestions $\alpha$ and calculate payment schemes $\pi$ that both depend on the exact realization $c$ of the demand of the truck drivers as can be seen from the optimization problems (\ref{optim2}) and (\ref{approx_optimization}), respectively, while the corresponding quantities under the CPURR scheme do not depend on $c$ as can be seen from (\ref{cpurr}). In other words, under the OPS and the AOPS schemes, the coordinator provides routing suggestions having complete information of the traffic demand. On the other hand, under the CPURR scheme, the truck drivers make individual routing decisions having incomplete information of the traffic demand. 

Based on the simulation results presented in Figure~\ref{deterministic_vs_stochastic}, it can be observed that even in the case where the truck drivers know the exact realization of the demand for the rest of the truck drivers, the expected total monetary cost and the total objective value are lower in the case where OPS and AOPS are applied rather than in the case where the CPURR scheme is applied. These results demonstrate that OPS and AOPS can be efficiently used in both a deterministic and a stochastic demand scenario.

\subsection{Computational Time}

In this subsection, we gradually increase the number of OD pairs that the truck drivers follow in the Sioux Falls network and we measure the computational time needed to solve the UE, the SO, the OPS, the AOPS and the CPURR problems. For the experiments of this subsection, we assumed that the cost of each link as well as the number of passenger vehicles at each link of the network are identical to the ones used in Section 4.1. We further assumed that the truck drivers choose their class among two available options with VOT $s_1 = 200 \frac{\$}{hr}$ and $s_2 = 50 \frac{\$}{hr}$, respectively. The weighting factors were chosen to be $\lambda = \mu = 0.9$. The demand of the truck drivers was assumed to be stochastic and take one of the two equiprobable values, namely $d_1$ and $d_2$\footnote{The values of $d_1$ and $d_2$ can be found in this \href{https://www.dropbox.com/s/xeqpt6aq2a7enes/OD_pairs_and_Truck_demands_TR_Part_C.txt?dl=0}{link}.}.

\begin{table*}[h]
\begin{center}
\begin{adjustbox}{max width=\textwidth}
\begin{tabular}{cc|c|c|c|c|c|}
\\
\cline{3-7} 
OD pairs&Metric&\textbf{UE}&\textbf{SO}&\textbf{OPS}&\textbf{AOPS}&\textbf{CPURR}\\
\hline
\multirow{4}{*}{{{4}}}&$E[T_{tr}]$&10545.7&10505.2&10505.5&10504.8&10520.7\\
&$E[T_{tr}^{mon}]$&23435.6&23105.2&23107.0&23139.3&23365.7\\
&$E[T_{s}]$&261346.5&261446.1&261444.2&261429.7&261337.4\\
&$E[T_{p}]$&250800.8&250940.9&250938.7&250924.9&250816.7\\
\hline
\multirow{4}{*}{{{8}}}&$E[T_{tr}]$&29463.3&28798.1&28801.9&28925.9&29111.5\\
&$E[T_{tr}^{mon}]$&62499.2&60282.2&60264.3&60459.3&61626.4\\
&$E[T_{s}]$&393917.8&393062.0&393055.2&392995.1&392876.2\\
&$E[T_{p}]$&364454.5&364263.9&364253.3&364069.1&363764.7\\
\hline
\multirow{4}{*}{{{12}}}&$E[T_{tr}]$&56222.9&54732.2&54719.0&55445.5&55392.5\\
&$E[T_{tr}^{mon}]$&114578.8&110117.8&110343.3&113907.9&112265.4\\
&$E[T_{s}]$&554185.9&549492.6&549476.0&551882.0&551676.9\\
&$E[T_{p}]$&497963.0&494760.4&494757.0&496436.5&496284.4\\
\hline
\multirow{4}{*}{{{16}}}&$E[T_{tr}]$&60264.6&58852.1&58847.2&59443.1&C.I.\\
&$E[T_{tr}^{mon}]$&122731.0&118475.3&118599.2&121315.8&C.I.\\
&$E[T_{s}]$&575701.4&571333.9&571307.1&573482.9&C.I.\\
&$E[T_{p}]$&515436.8&512481.8&512459.9&514039.8&C.I.\\
\hline
\multirow{4}{*}{{{20}}}&$E[T_{tr}]$&84675.5&82863.7&C.I.&83715.5&C.I.\\
&$E[T_{tr}^{mon}]$&169191.0&163705.0&C.I.&168081.6&C.I.\\
&$E[T_{s}]$&701513.8&696210.3&C.I.&699002.9&C.I.\\
&$E[T_{p}]$&616838.3&613346.6&C.I.&615287.5&C.I.\\
\hline
\end{tabular}
\end{adjustbox}
\end{center}
\caption{\label{OD_pairs}The expected total truck travel time $E[T_{tr}]$, the expected total truck monetary cost $E[T_{tr}^{mon}]$, the expected total network time $E[T_{s}]$ and the expected travel time of the passenger vehicles $E[T_{p}]$ of the UE, SO, OPS, AOPS and CPURR for different number of OD pairs, considering the 10 least congested routes per OD pair. The ratio of trucks in the network is shown in Table~\ref{ratio_of_trucks}.}
\end{table*}

\begin{table*}[t]
\begin{center}
\begin{adjustbox}{max width=\textwidth}
\begin{tabular}{c|c|c|c|c|c|}
\cline{2-6} 
OD pairs&\textbf{UE}&\textbf{SO}&\textbf{OPS}&\textbf{AOPS}&\textbf{CPURR}\\
\hline
\multirow{1}{*}{{{4}}}&6.9&8.9&32.5&10.0&33.9\\
\hline
\multirow{1}{*}{{{8}}}&30.7&24.9&337.4&59.1&1362.7\\
\hline
\multirow{1}{*}{{{12}}}&81.3&55.4&1389.0&128.2&6311.0\\
\hline
\multirow{1}{*}{{{16}}}&151.2&105.7&12163.5&179.4&C.I.\\
\hline
\multirow{1}{*}{{{20}}}&278.3&169.8&C.I.&446.4&C.I.\\
\hline
\end{tabular}
\end{adjustbox}
\end{center}
\caption{\label{Comp_Time}The computational time (in seconds) of the simulation experiments presented in Table~\ref{OD_pairs}.}
\end{table*}

In Table~\ref{OD_pairs}, we show the expected total truck travel time $E[T_{tr}]$, the expected total truck monetary cost $E[T_{tr}^{mon}]$, the expected total travel time of the network $E[T_{s}]$ and the expected total travel time of the passenger vehicles $E[T_{p}]$ of the UE, SO, OPS, AOPS and CPURR for different number of OD pairs, considering the 10 least congested routes per OD pair. In Table~\ref{Comp_Time}, we show the corresponding computational times. As can be seen from the results of Table~\ref{OD_pairs}, OPS always provides a solution close to the SO while AOPS outperforms both the CPURR scheme and the UE solution. Additionally, we observe that the CPURR scheme is the slowest approach since it becomes Computationally Intractable (C.I.) when 16 OD pairs are considered, while the OPS becomes C.I. when the truck drivers are assumed to follow 20 OD pairs. As far as it concerns the AOPS, it can be observed that even though it is not as efficient as the OPS in the evaluation metrics used, it remains computationally for a larger number of OD pairs. This result was expected since OPS calculates both the optimum routing strategy $\alpha^*_{OPS}$ and the optimum pricing scheme $\pi^*_{OPS}$ through the solution of the optimization problem (\ref{optim2}). On the other hand, under the AOPS, in the optimization problem (\ref{approx_optimization}), we only need to find the optimum routing strategy $\alpha^*_{AOPS}$ since the approximately optimum pricing scheme $\pi^{AOPS}_{c,j,w,r}$ can be easily calculated from the algebraic equation (\ref{payments}). This approximation makes AOPS significantly faster compared to OPS as can be seen from the results of Table~\ref{Comp_Time}. However, as the number of OD pairs increases, the computational time increases as well and after a certain number of OD pairs, even AOPS will become computationally intractable. Therefore, for realistic transportation networks with hundreds or thousands of OD pairs, more computationally efficient solutions need to be studied. 

We next show that the effect of the proposed methods on passenger vehicles is minimal, whereas for the trucks, it can be impactful. First, as can be seen in Table~\ref{ratio_of_trucks}, the number of trucks in the network is relatively small compared to the number of passenger vehicles. Additionally, as can be observed from the results of Table~\ref{OD_pairs}, the effect of truck routing schemes in the expected total travel time of the passenger vehicles is insignificant. More specifically, in the case of 4 OD pairs, OPS increases the expected travel time of the passenger vehicles in the network by $0.05\%$ compared to the UE. In the rest of the cases, both OPS and AOPS decrease the expected travel time of the passenger vehicles in the network. However, this reduction can be still considered small since it reaches up to $0.6\%$ in the case of 12 and 16 OD pairs. Therefore, we expect that the passenger vehicles will not react to such a minimal change of their traffic environment. On the other hand, the proposed pricing-and-routing schemes can offer significant benefits to trucks since they can reduce the expected total travel time of the truck drivers by $2.7\%$ and their corresponding expected total monetary cost by $3.7\%$ compared to the UE as can be observed from the results of Table~\ref{OD_pairs}.  

\begin{table*}[h]
\begin{center}
\begin{adjustbox}{max width=\textwidth}
\begin{tabular}{cc|c|}
\\
\cline{3-3} 
OD pairs&Realization&\textbf{Ratio (\%)}\\
\hline
\multirow{2}{*}{{{4}}}&$d_1$&2.82\\
&$d_2$&3.89\\
\hline
\multirow{2}{*}{{{8}}}&$d_1$&6.36\\
&$d_2$&6.96\\
\hline
\multirow{2}{*}{{{12}}}&$d_1$&8.64\\
&$d_2$&8.98\\
\hline
\multirow{2}{*}{{{16}}}&$d_1$&10.33\\
&$d_2$&10.40\\
\hline
\multirow{2}{*}{{{20}}}&$d_1$&11.88\\
&$d_2$&12.00\\
\hline
\end{tabular}
\end{adjustbox}
\end{center}
\caption{\label{ratio_of_trucks}The ratio of trucks in the network for the experimental results presented in Table~\ref{OD_pairs}.}
\end{table*}

\section{Conclusion.}
In this paper, we designed pricing-and-routing schemes that can be applied in a general transportation network to alleviate traffic congestion. We particularly focused on the design of a coordination system for the truck drivers in the case of stochastic OD demand considering their heterogeneity in VOT. In contrast to previous efforts that proposed class-anonymous pricing schemes using the idea of Congestion Pricing with Uniform Revenue Refunding (CPURR), or class-specific pricing schemes by making assumptions about the distribution that the VOT of the drivers might follow, we designed personalized (VOT-based) pricing-and-routing schemes where the users directly report their VOT to a central coordinator. More specifically, assuming that the users are asked to declare their OD pair and additionally pick their VOT from a set of $N$ available options, we proved the existence of Pareto-improving and revenue-neutral (budget balanced) on average pricing schemes that can additionally guarantee that every truck driver will have an incentive to truthfully declare his/her VOT. We showed that the Optimum Pricing Scheme (OPS) can be calculated by solving a nonconvex optimization problem and we additionally proposed an Approximately Optimum Pricing Scheme (AOPS) to approximate the solution of the OPS and reduce the computational time. Finally, we experimentally showed that both OPS and AOPS can significantly reduce the expected total travel time and the expected total monetary cost of the users and approximate the SO solution, while concurrently outperforming both the UE and the CPURR scheme, demonstrating the efficiency of VOT-based pricing schemes compared to class-anonymous pricing schemes.

There are several possible extensions of this work. First, both OPS and AOPS are route-based pricing schemes that cannot be easily implemented for passenger traffic. Therefore, personalized (VOT-based) pricing schemes which satisfy the properties of OPS and AOPS and can be applied to passenger traffic and multimodal transportation networks need to be investigated. Second, even though AOPS is more computationally efficient compared to OPS and CPURR, it is expected that as the number of OD pairs increases, AOPS will also become computationally intractable. Therefore, for realistic transportation networks with hundreds or thousands of OD pairs, more computationally efficient solutions need to be studied. Nevertheless, we believe that both OPS and AOPS can form the foundation and initial step in the design of VOT-based pricing schemes where the drivers directly report their VOT to a central authority that can also deal with the computational complexity of realistic transportation networks. Finally, the extension of the current work to the case where the planning horizon consists of multiple time windows and the drivers are asked to pick both their time window as well as their VOT, or to the case where it can be applied in real-time are of major importance.

\ACKNOWLEDGMENT{Funding: This work has been supported by the National Science Foundation Awards CPS \#1545130 and CNS-1932615.}


%

%
%
\begin{APPENDICES}

\newpage

\section{Notation Table}
\begin{table}[htbp]
\begin{center}
\begin{tabular}{|c|p{105mm}|}
\hline
Variable & Meaning \\\hline
$G$ & The transportation network as a graph\\
$V$ & Set of nodes in the network\\
$L$ & Set of links in the network\\
$R_j$ & Set of available routes connecting OD pair $j$\\
$m$ & Number of road segments in the network\\
$d_{j,w}^{c}$ & Demand of truck drivers belonging to the class $w$ with desired OD pair $j$ during the demand realization $c$\\
$\alpha_{w,r}^{c,j}$ & Proportion of truck drivers belonging to class $w$ with desired OD pair $j$ who follow route $r$ during the demand realization $c$\\
$s_w$ & Value of Time (VOT) of a truck driver belonging to class $w$\\
$X_{lp}$ & Number of passenger vehicles in the road segment $l$\\
$X_{lT}$ & Number of trucks in road segment $l$\\
$C_{lT}$ & Travel time of a truck driver traversing road segment $l$\\
$F^{w}_{j,r}$ & Expected travel time of a truck driver of class $w$ with desired OD pair $j$ who follows route $r$\\
$E[T_{tr}]$ & Expected total travel time of the truck drivers in the network\\
$E[T^{mon}_{tr}]$ & Expected total monetary cost of the truck drivers\\
$p_c$ & Probability of the demand realization $c$\\
$\alpha_{j,w,r}^{UE}$ & Proportion of truck drivers belonging to class $w$ with a desired OD pair $j$ who follow route $r$ at the UE\\
$J_{c,j,w,r}^{UE}$ & Travel time of a truck driver belonging to class $w$ with an OD pair $j$ who follows route $r$ during the demand realization $c$ at the UE\\
$E[T_p]$ & Expected total travel time of the passenger vehicles\\
$E[T_s]$ & Expected total travel time of the network \\
$\lambda$ & Weighting factor of the objective function \\
$\mu$ & Weighting factor of the objective function \\
$N$ & Number of classes with different VOT\\
$J^{M,c,j}_{w,r}$ & Travel time of a truck driver belonging to class $w$ with OD pair $j$ who follows route $r$ during the demand realization $c$ under the mechanism suggestions $M$\\
$\pi^{c,j}_{w,r}$ & Payments of truck drivers of class $w$ with desired OD pair $j$ who follow route $r$ during the demand realization $c$\\
$A^{UE}_{c,j}$ & Average travel time of a truck driver with OD pair $j$ during the demand realization $c$ at the UE\\
$F^{CP}_{j,w,r}$ & Expected total cost (travel time + payments expressed in time units) of a truck driver belonging to class $w$ with desired OD pair $j$ who follows route $r$ under the CPURR scheme\\
$J^{CP}_{c,j,w,r}$ & Travel time of a truck driver belonging to class $w$ with OD pair $j$ who follows route $r$ during the demand realization $c$ under the CPURR scheme\\
\hline
\end{tabular}
\end{center}
\caption{\label{notation_table}Notation used.}
\end{table}

\section{Discussion on the the Pareto-Improvement and the Truthfulness Properties\label{expost_analysis}}

In Sections 3.1 and 3.2, we designed the Optimum Pricing Scheme (OPS) and the Approximately Optimum Pricing Scheme (AOPS) and we guaranteed the Pareto-improvement, the truthfulness and the budget balance on average properties.

As also mentioned in Sections 3.1 and 3.2, using the formulations of OPS and AOPS, it is theoretically possible that for some demand realizations, some drivers are given a route with a higher total travel time (travel time + payments expressed in time units) compared to the average travel time at the UE. However, due to the fact that individual drivers only know the probability distribution of the demand and not the exact realization of it and hence they have incomplete information of the traffic conditions, at the time they make their decision, their expected cost under the mechanism suggestions $M$ is going to be lower than their corresponding cost at the UE and therefore, they will always be willing to participate in both coordinated schemes. At this point, we should note that it is possible to design an even more robust version for the Pareto-improvement property that guarantees that for every demand realization $c$ and for every route $r$, the total travel time of each driver (travel time + payments expressed in time units) is going to be lower under the mechanism suggestions compared his/her average travel time at the UE. We call this property {\it robust Pareto-improvement}.

As far as it concerns the truthfulness property, both OPS and AOPS guaranteed that every driver will have an incentive to truthfully declare his/her VOT since at the time they make their decision, their expected cost in the case they are truthful is lower than their corresponding cost in the case where they declared a different VOT than their actual one. We call this property {\it ex-ante truthfulness}. However, note that we are able to guarantee an even stronger version of the truthfulness property by ensuring that at each demand realization $c$, the average total travel time (travel time + payments expressed in time units) of a truck driver in case they are truthful is lower than their corresponding total travel time in the case where they declared a different VOT than their actual one. We call this property {\it ex-post truthfulness}.

Based on the aforementioned, we want to design a mechanism with the following properties.
\begin{itemize}
    \item A mechanism is {\it robust Pareto-improving} if:
    \begin{equation}
    \label{append_pareto}
        J^{M,c,j}_{w,r} + \frac{1}{s_w} \pi^{c,j}_{w,r} \leq A^{UE}_{c,j}, \forall c,j,w,r
    \end{equation}
    
    \item A mechanism is {\it ex-post truthful} if:
    \begin{equation}
    \label{append_truthful}
        \sum_{r \in R_j} \alpha^{c,j}_{i,r} (J^{M,c,j}_{i,r} + \frac{1}{s_i} \pi^{c,j}_{i,r}) \leq \sum_{r \in R_j} \alpha^{c,j}_{k,r} (J^{M,c,j}_{k,r} + \frac{1}{s_i} \pi^{c,j}_{k,r}), \; \forall c,j,i,k \\
    \end{equation}
    
    \item A mechanism is {\it budget balanced} on average if:
    \begin{equation}
    \label{append_BB}
        \sum_{c} \sum_{j=1}^v \sum_{w=1}^N \sum_{r \in R_j} p_c d^{w}_{c,j} \alpha^{c,j}_{w,r} \pi^{c,j}_{w,r} = 0
    \end{equation}
\end{itemize}

For convenience of the reader, we state that the payment scheme that we used in the AOPS formulation is given by:

\begin{equation}
\label{append_payments}
\begin{aligned}
\pi^{AOPS}_{c,j,w,r} = \ &s_{w} (A^{UE}_{c,j}- J^{M,c,j}_{w,r}) + \frac{s_{w}}{\sum_{l=1}^N s_l} \frac{E\big[T^{mon,M}_{tr}\big] - E\big[T^{mon,UE}_{tr}\big]}{\sum_{j=1}^v d^{w}_{c,j}}
\end{aligned}
\end{equation}

Now, let us formulate the following optimization problem:

\begin{equation}
\label{append_approx_optimization}
\begin{aligned}
& \underset{\alpha(\cdot)}{\text{minimize}}
& & \lambda (\mu E[T_{tr}(\alpha)] + (1-\mu)E[T_{p}(\alpha)])+(1-\lambda)E[T^{mon}_{tr}(\alpha)] \\
& \text{subject to}
& & E[T^{mon,M}_{tr}] \leq E[T^{mon,UE}_{tr}]\\
&&& Q_{i,k}^{c,j}(\alpha) \leq U_{i,k}^{c,j}(\alpha), \; \forall c,j,i,k \\
&&& \sum_{r \in R_j} \alpha^{c,j}_{w,r}=1, \; \forall c,j,w \\
&&& \alpha^{c,j}_{w,r}\geq 0, \; \forall c,j,w,r
\end{aligned}
\end{equation}
where $Q_{i,k}^{c,j}(\alpha)$ and $U_{i,k}^{c,j}(\alpha)$ are given by the following equations:
\begin{equation}
\label{append_cond2}
\begin{aligned}
Q_{i,k}^{c,j}(\alpha) = \bigg(1-\frac{s_k}{s_i}\bigg) A^{UE}_{c,j} + \frac{1}{\sum_{l=1}^N s_l} \frac{E[T^{mon,M}_{tr}] - E\big[T^{mon,UE}_{tr}\big]}{\sum_{j=1}^v d^i_{c,j}} 
\end{aligned}
\end{equation}
\begin{equation}
\label{append_cond3}
\begin{aligned}
U_{i,k}^{c,j}(\alpha) = \bigg(1-\frac{s_k}{s_i}\bigg) \sum_{r \in R_j} \alpha^{c,j}_{k,r} J^{M,c,j}_{k,r} + \frac{s_k}{s_i} \frac{1}{\sum_{l=1}^N s_l} \frac{E[T^{mon,M}_{tr}] - E\big[T^{mon,UE}_{tr}\big]}{\sum_{j=1}^v d^k_{c,j}} 
\end{aligned}
\end{equation}

Note that a solution to the optimization problem described by (\ref{append_approx_optimization})-(\ref{append_cond3}) always exists since the UE satisfies all of its constraints. Let $\alpha^{*}_{AOPS-EPT}$ be the optimum solution of the optimization problem described by (\ref{append_approx_optimization})-(\ref{append_cond3}).

In Proposition 1, we prove that the pricing-and-routing scheme described by (\ref{append_payments})-(\ref{append_cond3}) can guarantee the robust Pareto-improvement, the ex-post truthfulness and the budget balance on average properties. We call this scheme Approximately Optimum Pricing Scheme with Ex-Post Truthfulness guarantees (AOPS-EPT). Note that the only difference between AOPS-EPT and AOPS is at the second constraint of the corresponding optimization problems. More specifically, observe that the second constraint of (\ref{approx_optimization}) which corresponds to AOPS, is calculated by taking the expectation over the different demand realizations which makes the truthfulness property hold ex-ante, i.e. at the time they make their decision, the expected cost of the truck drivers in the case they are truthful is lower than their corresponding cost in the case where they declared a different VOT than their actual one. On the other hand, the second constraint of (\ref{append_approx_optimization}) which corresponds to AOPS-EPT, holds for every possible demand realization $c$, making the truthfulness property hold ex-post.

\begin{proposition}
The pair ($\alpha_{AOPS-EPT}^{*}$, $\pi^{AOPS}_{c,j,w,r}$) makes every truck driver better-off compared to the UE (robust Pareto-improvement), guarantees that every user will have an incentive (ex-post) to truthfully declare his/her VOT and leads to a budget balanced on average mechanism.
\end{proposition}
{{ \bf Proof.} Similar to the proof of Theorem 1, substituting (\ref{append_payments}) into (\ref{append_pareto}), we get:
\begin{equation*}
\begin{aligned} J^{M,c,j}_{w,r} + A^{UE}_{c,j} - J^{M,c,j}_{w,r} + \frac{1}{\sum_{l=1}^N s_l} \frac{E[T^{mon,M}_{tr}] - E\big[T^{mon,UE}_{tr}\big]}{\sum_{j=1}^v d^{w}_{c,j}} \leq A^{UE}_{c,j} \Leftrightarrow \\
 \Leftrightarrow \frac{1}{\sum_{l=1}^N s_l} \bigg(\frac{E[T^{mon,M}_{tr}] - E\big[T^{mon,UE}_{tr}\big]}{\sum_{j=1}^v d^{w}_{c,j}}\bigg) \leq 0
\end{aligned}
\end{equation*}
which holds true if and only if $E[T^{mon,M}_{tr}] \leq E[T^{mon,UE}_{tr}]$
which is equivalent to the first constraint of (\ref{append_approx_optimization}). Additionally, a user will have an incentive (ex-post) to truthfully declare his/her VOT if (\ref{append_truthful}) holds. Therefore, substituting (\ref{append_payments}) into (\ref{append_truthful}), we get:
\begin{equation*}
\begin{aligned}
\sum_{r \in R_j} \alpha^{c,j}_{i,r} \bigg(J^{M,c,j}_{i,r} + A^{UE}_{c,j} - J^{M,c,j}_{i,r} + \frac{1}{\sum_{l=1}^N s_l} \frac{E[T^{mon,M}_{tr}] - E\big[T^{mon,UE}_{tr}\big]}{\sum_{j=1}^v d^{i}_{c,j}} \bigg) \leq \\ \leq
\sum_{r \in R_j} \alpha^{c,j}_{k,r} \bigg(J^{M,c,j}_{k,r} + \frac{s_k}{s_i} \bigg( A^{UE}_{c,j} - J^{M,c,j}_{k,r} + \frac{1}{\sum_{l=1}^N s_l} \frac{E[T^{mon,M}_{tr}] - E\big[T^{mon,UE}_{tr}\big]}{\sum_{j=1}^v d^{k}_{c,j}} \bigg)\bigg) \Leftrightarrow \\
\Leftrightarrow  
\bigg(1-\frac{s_k}{s_i}\bigg) A^{UE}_{c,j} + \frac{1}{\sum_{l=1}^N s_l} \frac{E[T^{mon,M}_{tr}] - E\big[T^{mon,UE}_{tr}\big]}{\sum_{j=1}^v d^i_{c,j}} \leq \\
\leq  \bigg(1-\frac{s_k}{s_i}\bigg) \sum_{r \in R_j} \alpha^{c,j}_{k,r} J^{M,c,j}_{k,r} + \frac{s_k}{s_i} \frac{1}{\sum_{l=1}^N s_l} \frac{E[T^{mon,M}_{tr}] - E\big[T^{mon,UE}_{tr}\big]}{\sum_{j=1}^v d^k_{c,j}} 
\end{aligned}
\end{equation*}
where the last inequality is equivalent to the second constraint of (\ref{append_approx_optimization}). Last, a mechanism is budget balanced on average if (\ref{append_BB}) holds. In Theorem 1, we proved that (\ref{append_payments}) satisfies the budget balance on average property and this concludes the proof. \hfill \qed}

There are three main reasons that justify our choice to use the AOPS instead of the AOPS-EPT in our analysis. First, since the drivers only know the probability distribution of the demand and not the exact realization of it, they have incomplete information of the traffic conditions. AOPS is still sufficient to guarantee that the drivers will have an incentive to participate in the mechanism and truthfully declare their VOT since at the time they make their decision, their expected cost in the case they are truthful is lower than their corresponding cost in the case where they declared a different VOT than their actual one. Second, using AOPS-EPT would be less computationally efficient since the second constraint of (\ref{append_approx_optimization}) should hold for every possible demand realization $c$ which significantly increases the number of constraints compared to the optimization problem (\ref{approx_optimization}). Third, one can expect that the optimum solution calculated by AOPS-EPT would be less efficient compared to the one calculated by AOPS since the size of the feasible region over which we optimize is smaller.

\end{APPENDICES}

\newpage


\bibliographystyle{ifacconf.bst} 
\bibliography{ifacconf.bib} 


\end{document}